\documentclass[12pt]{amsart}
\usepackage[dvips]{epsfig}
\usepackage{graphics}
\usepackage{latexsym}
\usepackage{verbatim}
\usepackage{amsmath}
\usepackage{amsthm}
\usepackage{amssymb}
\newtheorem{theorem}{Theorem}[section]
\newtheorem{lemma}[theorem]{Lemma}

\def\Remark{\medskip\noindent{\bf Remark: }}
\def\Remarks{\medskip\noindent{\bf Remarks: }}

\def\bul{{$\bullet$\hspace*{2mm}}}

\newcommand{\ens}[1]{\mathbb{#1}}
\newcommand{\ron}[1]{\mathcal{#1}}
\newcommand{\N}{\mathbb{N}}

\newcommand{\R}{\mathbb{R}}
\newcommand{\C}{\mathbb{C}}

\def\cal{\mathcal}

\def\derpar#1#2{\frac{\partial#1}{\partial#2}}
\def\var{\varepsilon}
\def\signcm{\bigskip\bigskip\hspace{80mm}
\vbox{{\sc C. Mouhot\par\vspace{3mm}
CEREMADE, Universit\'e Paris IX Dauphine \par
Place du Maréchal de Lattre de Tassigny \par
75775 Paris Cedex 16 \par
FRANCE\par\vspace{3mm}
e-mail:} cmouhot@ceremade.dauphine.fr }}

\begin{document}

\title[Coercivity estimates for the Boltzmann and Landau operators]
{Explicit coercivity estimates for the linearized Boltzmann and Landau operators}

\author{Cl\'ement Mouhot}

\maketitle

\hyphenation{bounda-ry rea-so-na-ble be-ha-vior pro-per-ties
cha-rac-te-ris-tic}

\begin{abstract} We prove explicit coercivity estimates for 
the linearized Boltzmann and Landau operators, for a general 
class of interactions including any inverse-power law interactions, 
and hard spheres. The functional spaces of these  
coercivity estimates depend on the collision kernel of these operators. 
They cover the spectral gap estimates for the linearized Boltzmann 
operator with Maxwell molecules, improve these estimates 
for hard potentials, and are the first explicit coercivity estimates for soft 
potentials (including in particular the case of Coulombian interactions). 
We also prove a regularity property for 
the linearized Boltzmann operator with non locally 
integrable collision kernels, and we deduce from it a new proof of 
the compactness of its resolvent for hard potentials without angular cutoff. 
\end{abstract}
\bigskip



\tableofcontents

\section{Introduction and main results}
\setcounter{equation}{0}

This paper is devoted to the study of the linearized 
Boltzmann and Landau collision operators. In this work we 
shall obtain new quantitative coercivity estimates for these operators. 
Before we explain our methods and results in more details, let us introduce 
the problem in a precise way. 

\subsection{The models}
The {\em Boltzmann equation} describes the behavior of a dilute 
gas when the only interactions taken into account are binary elastic collisions. 
It reads in $\R^N$ ($N \ge 2$) 
  \begin{equation*}
  \derpar{f}{t} + v \cdot \nabla_x f = Q^{{\ron B}}(f,f),
  \end{equation*}
where $f(t,x,v)$ stands for the time-dependent probability density of particles 
in the phase space. The $N$-dimensional {\em Boltzmann operator} $Q^{{\ron B}}$ 
is a quadratic operator, which is local in $(t,x)$. 
The time and position are only parameters and therefore 
shall be omitted in the sequel: the functional 
estimates proved in this paper are all local in $(t,x)$. 
This operator acts on $f(v)$ by 
  \begin{equation*}
  Q^{{\ron B}} (f,f)(v) = \int_{\R^N \times \ens{S}^{N-1}}  B(|v-v_*|, \cos \theta) \, 
  \big[ f'_* f' - f_* f \big] \, dv_* \, d\sigma
  \end{equation*} 
where we have used the shorthand $f = f(v)$, $f_* = f(v_*)$, $f ^{'} = f(v')$, 
$f_* ^{'} = f(v_* ^{'})$. In this formula, $v'$, $v' _*$ and $v$, $v_*$ are 
the velocities of a pair of particles before and after collision, they are related by 
  \begin{equation*}
  v' = (v+v_*)/2 + (|v-v_*|/2) \sigma, \qquad
  v'_* = (v+v_*)/2 - (|v-v_*|/2) \sigma.
  \end{equation*}
The {\em collision kernel} $B$ is a non-negative function which 
only depends on the {\em relative velocity} $|v-v_*|$ and 
the {\em deviation angle} $\theta$ through 
$\cos \theta = k \cdot \sigma$ where $k = (v-v_*)/|v-v_*|$. 

In the case of long-distance interactions, collisions occur 
mostly for very small $\theta$. When all collisions become 
concentrated on $\theta =0$, one obtains by the so-called 
{\em grazing collision limit} asymptotic 
(see for instance \cite{ArBu:91,DeLu:92,Desv:asBE:92,Gou:97,Vill:nocu:98,AlVi:04}) 
the {\em Landau operator} 
  \begin{equation*}
  Q^{{\ron L}} (f,f)(v) = \nabla _v \cdot \left( \int_{\R^N} 
  {\bf A}(v-v_*) \big[ f_* \left( \nabla f \right) 
  - f \left( \nabla f \right)_* \big] \, dv_* \right),
  \end{equation*} 
with ${\bf A} (z) = |z|^2 \, \Phi(|z|) \, {\bf P}(z)$, $\Phi$ is a non-negative function, 
and ${\bf P}(z)$ is the orthogonal projection onto $z^\bot$, {\it i.e.},  
 \begin{equation*}
 \left( {\bf P}(z) \right) _{i,j} = \delta_{i,j} - \frac{z_i z_j}{|z|^2}.
 \end{equation*}
This operator is used for instance in models of plasma in the case of a Coulomb 
potential where $\Phi(|z|) = |z|^{-3}$ in dimension $3$ 
(for more details see~\cite[Chapter~1, Section~1.7]{Vi:hand} and the references therein). 
Indeed for Coulombian interactions the Boltzmann collision operator 
does not make sense anymore (see~\cite[Annex~I, Appendix]{Vi:habil}). 

Boltzmann and Landau collision operators have the fundamental properties of
preserving mass, momentum and energy ($Q$ denotes $Q^{{\ron B}}$ or 
$Q^{{\ron L}}$) 
  \begin{equation*}
  \int_{\R^N} Q (f,f)\, \phi(v)\,dv = 0, \quad
  \phi(v)=1,v,|v|^2. 
  \end{equation*}
Moreover they satisfy well-known Boltzmann's $H$ theorem, which writes formally 
  \begin{equation*} 
  - \frac{d}{dt} \int_{\R^N} f \log f \, dv = - \int_{\R^N} Q(f,f)\log(f) \, dv \geq 0.
  \end{equation*}
The functional $- \int f \log f$ is the entropy of the solution. 
The $H$ theorem implies formally that any equilibrium
distribution, {\it i.e.}, any distribution which maximizes the entropy,
has the form of a locally Maxwellian distribution
  \begin{equation*}
  M(\rho,u,T)(v)=\frac{\rho}{(2\pi T)^{N/2}}
  \exp \left\{ - \frac{\vert u - v \vert^2} {2T} \right\}, 
  \end{equation*}
where $\rho,\,u,\,T$ are the mass, momentum and temperature of the gas
  \begin{equation*}
  \rho = \int_{\R^N}f(v)\, dv, \quad u =
  \frac{1}{\rho}\, \int_{\R^N}vf(v) \, dv, \quad T = {1\over{N\rho}}
  \, \int_{\R^N}\vert u - v \vert^2f(v) \, dv.
  \end{equation*}
For further details on the physical background and derivation of 
the Boltzmann and Landau equations we refer to~\cite{Ce88,CIP94,LandLif,Vi:hand}.

\subsection{Linearization}
Consider the linearization process $f = M (1+h)$ around the Maxwellian 
equilibrium state denoted by $M$. It yields the {\it linearized Boltzmann operator} 
  \begin{equation*}
  L^{{\ron B}} h(v)=\int_{\R^N \times \ens{S}^{N-1}} B(|v-v_*|,\cos \theta) \, 
  M(v_{*}) \left[ h^{'}_{*}+h^{'} - h_{*} - h \right] \, dv_{*} \, d\sigma,
  \end{equation*}
with (up to a normalization and without restriction) $M(v) = e^{-|v|^2}$. The normalization corresponds 
to the choice of mass $\pi^{N/2}$, momentum $0$ and temperature $1/2$. 
For the sake of simplicity we shall give all the statements and proofs under the 
normalization. It is explained briefly in Subsection~\ref{subsec:norm} how to obtain 
the dependency of the estimates we shall establish in terms of the mass, 
momentum and energy of the equilibrium in the importance case where the collision kernel 
depends on the relative kernel in a polynomial way. 

The operator $L^{{\ron B}}$ is self-adjoint on the Hilbert space $L^2(M)$, which is defined by 
the general convention ($W$ is any positive weight function)
  \[ L^2(W) = \left\{ h: \R^N \rightarrow \R \mbox{ measurable s. t. } 
                      \|h\|_{L^2(W)} ^2 := 
                      \int_{\R^N} h(v)^2 \, W(v) \, dv <+\infty \right\}. \]
Since we shall need it in the sequel, we define similarly 
  \begin{multline*}
  H^1(W) = \Bigg\{ h: \R^N \rightarrow \R \mbox{ measurable s. t. } 
                      \|h\|_{H^1(W)} ^2 := 
                      \int_{\R^N} h(v)^2 \, W(v) \, dv \\
                      + \int_{\R^N} |\nabla h(v)|^2 \, W(v) \, dv<+\infty \Bigg\}. 
  \end{multline*}

The Dirichlet form of $-L^{{\ron B}}$ in this space satisfies 
  \begin{eqnarray*}
  D^{{\ron B}} (h) &=& - \left( h,L^{{\ron B}} h \right)_{L^2(M)}\\ 
  &=& \frac{1}{4}\int_{\R^N \times \R^N \times \ens{S}^{N-1}} 
  B(|v-v_*|, \cos \theta) \, \left[ h^{'}_{*}+h^{'}-h_{*}-h \right]^2 
  M\, M_{*} \, dv \, dv_* \, d\sigma.
  \end{eqnarray*}
It is non-negative, which implies that the spectrum 
of $L^{{\ron B}}$ in $L^2(M)$ is included in $\R_-$. 
As it shall be useful in the sequel for the study of $L^{{\ron B}}$, 
we define the {\em collision frequency} (in $[0,+\infty]$) by 
 \[ \nu (v) = \int_{\R^N \times \ens{S}^{N-1}} B(|v-v_*|,\cos \theta) \, M(v_*) \, dv_*. \]

The same linearization process yields the {\it linearized Landau operator} 
  \begin{equation*}
  L^{{\ron L}} h(v)=M(v) ^{-1} \, \nabla_v \cdot \left( \int_{v_* \in \R^N} 
  {\bf A}(v-v_*) \big[ \left( \nabla h \right)  - \left( \nabla h \right)_* \big] 
  M M_* \, dv_* \right),
  \end{equation*}
which is self-adjoint on $L^2(M)$, and whose Dirichlet form satisfies 
  \begin{eqnarray*}
  D^{{\ron L}} (h)& =& - \left( h,L^{{\ron L}}h \right)_{L^2(M)} \\
  &=&\frac{1}{2}\int_{\R^N}\int_{\R^N}\Phi(v-v_*)|v-v_*|^2
  \big\|{\bf P}(v-v_*) \big[ \left( \nabla h \right) - \left( \nabla h \right)_*\big] \big\|^2 
  M\, M_{*} \, dv_{*}\, dv.
  \end{eqnarray*}
It is also non-negative, which implies that the spectrum of $L^{{\ron L}}$ in $L^2(M)$ is 
included in $\R_-$. 
The null space of the two operators $L^{{\ron L}}$ and $L^{{\ron B}}$ is 
  \[ N(L^{{\ron B}}) = N(L^{{\ron L}}) = \mbox{Span}\{1,v_1, \dots, v_N,|v|^2\} \]
(note that it is independent on the collision kernel). 
These two properties -- the fact that the time-derivative of the $L^2(M)$ norm 
is negative and the fact that the only functions which cancel this derivative are 
the collision invariants -- correspond to the $H$ theorem at the linearized level.  
We denote by $\Pi$ the orthogonal projection on this null space in $L^2(M)$. 

\subsection{Assumptions on the collision kernel} 
 \begin{itemize}
 \item In the case of the linearized Boltzmann operator, $B$ takes the product form
   \begin{equation}\label{eq:hyptens}
   B(|v-v_*|, \cos \theta) = \Phi(|v-v_*|) \, b(\cos \theta),
   \end{equation}
 where $\Phi$ and $b$ are non-negative functions. 
 \item The kinetic part $\Phi$ is bounded by a power-law:
   \begin{equation}\label{eq:hypPhi}
   \forall \, r \ge 0, \ \ \  C_\Phi \, r^\gamma \le \Phi(r) \le C' _\Phi \, r^\gamma.
   \end{equation}
 where $\gamma \in (-N,1]$ for the linearized Boltzmann operator, 
 or $\gamma \in [-N,1]$ for the linearized Landau operator, and 
 $C_\Phi, C' _\Phi >0$ are some constants. 
 \item In the case of the linearized Boltzmann operator, the angular part $b$ satisfies
   \begin{equation}\label{eq:hypb}
   C_b = \inf_{\sigma_1, \sigma_2 \in \ens{S}^{N-1}} 
   \int_{\sigma_3 \in \ens{S}^{N-1}} \min \{  b(\sigma_1 \cdot \sigma_3), 
   b(\sigma_2 \cdot \sigma_3) \} \, d\sigma_3 > 0. 
   \end{equation}
This quite unusual (but satisfied for 
all physical models) assumption was defined and used in the work~\cite{BaMo} in order to 
obtain explicit spectral gap estimates for hard potentials.   
 \item In the particular case of the linearized Boltzmann operator with a 
 non locally integrable collision kernel, we shall assume the following 
 control in order to derive coercivity estimates in Sobolev spaces
   \begin{equation}\label{eq:hypbnc}
   \forall \, \theta \in (0,\pi], \ \ \  \frac{c_b}{\theta^{N-1+\alpha}} \le 
                         b(\cos \theta) \le \frac{c_b '}{\theta^{N-1+\alpha}}
   \end{equation}
 for some constants $c_b, c' _b >0$ and $\alpha \in [0,2)$ (note that  
 assumption~\eqref{eq:hypbnc} implies straightforwardly assumption~\eqref{eq:hypb}). 
 The goal of this control is to measure the strength of the angular singularity, which is related 
 to the regularity properties of the collision operator (see \cite{ADVW} for instance). 
 \end{itemize}

For the linearized Boltzmann operator, our assumptions 
cover in dimension $3$ the hard spheres collision kernel 
$B(|v-v_*|, \cos \theta)= \mbox{cst} \, |v-v_*|$ (which 
satisfies~(\ref{eq:hyptens}-\ref{eq:hypPhi}-\ref{eq:hypb})), and collision kernels  
deriving from interaction potentials behaving like inverse-power laws 
(which satisfy~(\ref{eq:hyptens}-\ref{eq:hypPhi}-\ref{eq:hypbnc})). More precisely 
for an interaction potential $V(r) = \mbox{cst} \, r^{-(s-1)}$, $B$ satisfies the 
assumptions~(\ref{eq:hyptens}-\ref{eq:hypPhi}-\ref{eq:hypbnc}) 
with the formulas $\gamma = (s-5)/(s-1)$ and $\alpha = 2/(s-1)$ (see~\cite{Ce88}). 
One traditionally denotes by {\em hard potentials} the case $s > 5$ 
(for which $0< \gamma <1$), {\em Maxwell molecules} the case $s=5$ (for which 
$\gamma=0$), and {\em soft potentials} the case $2<s<5$ (for which $-N<\gamma<0$). 
In the case of an {\em angular cutoff}, $b$ is assumed artificially to 
be integrable.  

For the linearized Landau operator, our assumptions are satisfied for any 
collision kernel $\Phi(|v-v_*|) = |v-v_*|^\gamma$, $-N \le \gamma \le 1$, including 
in particular in dimension $3$ the most important case of the 
linearized Landau-Coulomb operator when $\gamma=-3$ (corresponding 
to Coulombian interactions between charged particles). By analogy with 
the linearized Boltzmann operator, one traditionally denotes by 
{\em hard potentials} the case $\gamma >0$, {\em Maxwell 
molecules} the case $\gamma=0$, and {\em soft potentials} the case $-N \le \gamma<0$. 

\Remark The assumption~\eqref{eq:hyptens} is made for a sake of 
simplicity. Indeed, one could easily adapt 
the proofs in Section~\ref{sec:Bolt} to relax this assumption to 
upper and lower bounds on $B$ of the multiplicative type~\eqref{eq:hyptens}. 
For more general $B$ without decoupling between the arguments, 
more technical conditions would be needed to apply our strategy. 
\medskip

Let us comment on the domains ${\ron D}\left(L^{{\ron B}}\right)$ and 
${\ron D}\left(L^{{\ron L}}\right)$ of the linearized Boltzmann and Landau operators 
according to the previous assumptions: 

\bul In the case of the linearized Boltzmann collision operators with locally integrable collision kernel, 
it is straightforward that the domain is $L^2(\langle v \rangle^{\gamma^+} M)$ where 
$\gamma^+$ is the nonnegative part of $\gamma$. Indeed as discussed in Section~\ref{sec:Bolt} 
it is well-known (see~\cite{Grad63} for instance) that $L^{{\ron B}}$ then splits between a bounded 
part and a multiplicative part, which is the multiplication operator by the collision frequency. 

\bul In the case of the linearized Landau operator, obvious estimates show that the domain 
contains $H^1(\langle v \rangle ^{\gamma+2} M) \cap L^2(M)$. 
The maximal domain is slightly bigger than this space, it is obtained in~\cite{Guo:Land}:
$$ 
h \in L^2(M) \mbox{ s.t. } \| h \|_{L^2(\langle v \rangle^{\gamma+2} M)} + 
\left\| {\bf P}(v) \nabla_v h \right\|_{L^2(\langle v \rangle^{\gamma+2} M)} + 
\left\| (\mbox{Id} - {\bf P}(v)) \nabla_v h \right\|_{L^2(\langle v \rangle^\gamma M)} < + \infty
$$
where ${\bf P}(v)$ is the orthogonal projection on $v^\bot$ defined above. 

\bul In the case of the linearized Boltzmann operator with non locally integrable collision kernel, 
the question of characterizing the domain in terms of some 
known functional spaces is still open. However one can prove easily (using a Taylor 
expansion and the kind of changes of variables of Subsection~\ref{subsec:nc}) that 
$$ 
\| L^{{\ron B}} h \|_{L^2(M)} \le C \, \| h \|_{H^1\left(\langle v \rangle^{(\gamma+2)^+} M\right)}
$$
which shows that the domain includes $H^1(\langle v \rangle^{(\gamma+2)^+} M)$. 

\subsection{Motivation}

We refer to~\cite{BaMo} and the references therein for a discussion about the interest of 
spectral gap estimates for the linearized Boltzmann and Landau operators and 
some review. Let us just recall that since Grad~\cite{Grad63}, 
spectral gap estimates is known to exist at least for Maxwell molecules, 
hard potentials (with or without angular cutoff) and hard spheres. 
However, apart for the case of Maxwell molecules, for 
which the linearized Boltzmann operator is diagonalized explicitly 
in~\cite{WCUh:LBE:70,Bo88}, the classical proof of the existence 
of a spectral gap by Grad is based on non-constructive arguments 
and does not provide any estimate. 
In~\cite{BaMo}, it is given a new method to obtain explicit spectral gap 
estimates for any $\gamma \ge 0$ (including the so-called 
hard potentials -- with or without angular cutoff -- and hard spheres cases). 
This method relies on a geometrical argument 
on the whole collision operator, with no need of splitting or angular cutoff 
assumptions. The result is also extended in the same work 
to the linearized Landau operator (for $\gamma \ge 0$) by a grazing collision asymptotic.

As for soft potentials in dimension $N=3$ with $-1 < \gamma <0$ and  
angular cutoff, it was proved in~\cite{Cafl:LBE1:80} that the Boltzmann linearized operator 
has no spectral gap. But if one allows a loss on the algebraic weight of the norm, 
it was proved in~\cite{GoPo:89} (in dimension $N=3$ with $-2 < \gamma <0$ and angular cutoff) 
a ``degenerated spectral gap'' result of the form 
  \begin{equation}\label{GP}
  D^{{\ron B}} (h) \ge C \, \big\| \big[h - \Pi(h) \big]\langle v \rangle^{\gamma/2} \big\|^2 _{L^2 (M)}  
  \end{equation} 
where we have denoted $\langle \cdot \rangle = \sqrt{1+|\cdot|^2}$. 
The proof was based on Weyl's Theorem about 
compact perturbation of the essential spectrum, and it leads to non explicit 
constants. The inequality~\eqref{GP} was then extended to the full range 
$-N < -\gamma <0$ in~\cite{Guo:soft} by a similar non-constructive approach. 
\smallskip

In this work we shall extend and complete the works~\cite{Cafl:LBE1:80,GoPo:89,Guo:soft} 
and \cite{BaMo} by 
  \begin{itemize}
  \item giving a constructive proof of~\eqref{GP} for soft potentials, with explicit estimate 
  on the constant (note that all through this paper the word ``explicit" refers to the fact that all steps 
  in our proofs are constructive and explicit computations of the constants could be extracted 
  from them);
  \item extending it to hard potentials ($\gamma>0$)  
  (note that for hard potentials this estimate is stronger than the usual spectral gap estimate); 
  \item extending this approach to the linearized Landau operator by proving 
  coercivity estimates in $H^1$ with a weight corresponding to the collision kernel; 
  \item giving a coercivity result in local Sobolev spaces for the linearized Boltzmann operator  
  with a non locally integrable collision kernel, and discussing the consequence on its spectrum.
  \end{itemize}

There are the first explicit coercivity estimates for weak collision interaction 
models such as the linearized Boltzmann for soft potentials or the linearized 
Landau-Coulomb operator (that is the linearized Landau operator 
in the case of Coulombian interactions). These estimates are crucial for instance for the 
construction of a quantitative perturbative theory near equilibrium and for 
obtaining quantitative rates of convergence to equilibrium, which 
shall be discussed in a forthcoming work. 

In addition, the regularity coercivity estimates derived 
for the linearized Boltzmann operator allow to give a new simpler proof of the 
compactness of the resolvent for hard potentials without angular cutoff 
(which implies that the spectrum is purely discrete in this case). 
This revisits partly results in~\cite{Pao,Klaus}.  
We recall that in~\cite{Pao}, it was proved that the linearized Boltzmann operator has 
compact resolvent for interaction potentials $V(r)=\mbox{cst} \, r^{-(s-1)}$ with $3<s<+\infty$ 
in dimension $N=3$. As a byproduct, our results give a new proof 
of this fact in the case $s \ge 5$ (hard potentials). 

\subsection{Statements of the results}

We now state our main results:
 \begin{theorem}[The linearized Boltzmann operator] \label{theo:Bolt}
 Under the assumptions~\eqref{eq:hyptens}, \eqref{eq:hypPhi}, \eqref{eq:hypb}, 
 the linearized Boltzmann operator $L^{{\ron B}}$ with collision kernel $B = \Phi \, b$ 
 satisfies
   \begin{equation}\label{coBo}
   \forall \, h \in {\ron D}\left(L^{{\ron B}}\right), \quad 
   D^{{\ron B}} (h) \ge C^{{\ron B}} _\gamma \, \big\| \big[h - \Pi(h)\big] \, 
                              \langle v \rangle ^{\gamma/2} \big\|_{L^2(M)} ^2,
   \end{equation}
 where $C^{{\ron B}} _\gamma$ is an explicit constant depending only on 
 $\gamma$, $C_\Phi$, $C_b$, and the dimension $N$. 
 \end{theorem}

\Remarks 

1. When the collision kernel is locally integrable, the collision 
frequency $\nu$ is finite, bounded from below, and asymptotically 
equivalent to $ \langle v \rangle ^\gamma$, 
and the estimate~\eqref{coBo} can be written in the following form 
  \[ D^{{\ron B}} (h) \ge \bar C^{{\ron B}} _\gamma \, \big\| \big[h - \Pi(h)\big] \, 
                              \nu^{1/2} \big\|_{L^2(M)} ^2 \]
for some explicit constant $\bar C^{{\ron B}} _\gamma >0$. 
\smallskip

2. When the collision kernel is not locally integrable and $b$ satisfies~\eqref{eq:hypbnc}, 
a natural conjecture would be that the estimate~\eqref{coBo} improves into 
  \[ D^{{\ron B}} (h) \ge C^{{\ron B}} _{\gamma,\alpha} \, \big\| \big[h - \Pi(h)\big] \, 
                              \langle v \rangle ^{\gamma/2} \big\|_{H^{\alpha/2}(M)} ^2, \]
where $\alpha \in [0,2)$ is the order of angular singularity, defined in \eqref{eq:hypbnc}, 
and $H^{\alpha/2} (M)$ is the Sobolev space defined by 
$H^{\alpha/2} (M) = \{ h \in L^2(M) \mbox{ s. t. } (1-\Delta)^{-\alpha/4}h \in L^2(M) \}$.  
We were not able to obtain this coercivity estimate, however we give in the following 
theorem its consequence in terms of local regularity. We denote by 
$H^{\alpha/2} _{\mbox{{\scriptsize loc}}}$ the space of functions whose restriction to any 
compact set $K$ of $\R^N$ belongs to 
$H^{\alpha/2}(K) = \{ h \in L^2(K) \mbox{ s. t. } (1-\Delta)^{-\alpha/4}h \in L^2(K) \}$ 
(here $L^2(K)$ denotes the space of functions square integrable on $K$).  
\medskip

 \begin{theorem}[The linearized Boltzmann operator for long-range interactions]\label{theo:Bolt:nc}
 Under the assumptions~\eqref{eq:hyptens}, \eqref{eq:hypPhi}, \eqref{eq:hypbnc}, 
 the linearized Boltzmann operator $L^{{\ron B}}$ with collision kernel $B = \Phi \, b$ 
 satisfies \eqref{coBo} and 
   \begin{equation}\label{coBonc}
   \forall \, h \in {\ron D}\left(L^{{\ron B}}\right), \quad 
   D^{{\ron B}} (h) \ge 
   C^{{\ron B}} _{\gamma,\alpha} \, \|h - \Pi(h) \|_{H^{\alpha/2} _{\mbox{{\scriptsize {\em loc}}}}} ^2,
   \end{equation}
 where $C^{{\ron B}} _{\gamma,\alpha}$ is an explicit constant depending only on 
 $\gamma$, $\alpha$, $C_\Phi$, $C_b$, $c_b$ and the dimension $N$. 
 \end{theorem} 

\Remarks 

1. When $\gamma>0$ and $\alpha>0$, one can deduce 
straightforwardly from Theorem~\ref{theo:Bolt:nc} that the operator $L^{{\ron B}}$ has 
compact resolvent, which implies that 
its spectrum is purely discrete in this case. 
Indeed let us pick any $\xi \in \C$ such that $L^{{\ron B}}-\xi$ is invertible (such $\xi$ 
exists since the operator is self-adjoint for instance), and let us denote by 
$R(\xi) = \big( L^{{\ron B}}- \xi \big)^{-1}$ the resolvent at this point. 
For any sequence $(g_n)_{n \ge0}$ bounded 
in $L^2(M)$, we can define the sequence $h_n = R(\xi)(g_n)$ which is also bounded in $L^2(M)$ 
since the operator $R(\xi)$ is bounded. We have:  
  \[ \forall \, n \ge 0, \ \ \ L^{{\ron B}}(h_n) = g_n + \xi \, h_n \]
and so the sequence $L^{{\ron B}}(h_n)$ is bounded in $L^2(M)$. It follows that 
the sequence $D^{{\ron B}}(h_n)$ is bounded in $\R$, and we deduce from the 
coercivity estimates above that 
the sequence $(h_n)_{n \ge 0}$ is bounded in 
$L^2(\langle v \rangle^\gamma M) \cap H^{\alpha/2} _{\mbox{{\scriptsize loc}}}$. 
When $\alpha >0$ and $\gamma > 0$, it implies that 
it has a cluster point in $L^2(M)$ by 
Rellich-Kondrachov compactness Theorem. Thus the operator $R(\xi)$ is 
compact. By classical arguments (see~\cite{Kato} for instance), it implies 
that the resolvent $R(\xi)$ is compact at every $\xi \in \C$ for which it is defined,  
and that the spectrum of $L^{{\ron B}}$ is purely discrete. 
\smallskip

2. The case $\alpha=0$ could probably 
be treated in the same spirit as in Theorem~\ref{theo:Bolt:nc}, 
using for the coercivity estimate a functional space 
controlling logarithmic derivatives defined by the norm $\|h \log(1-\Delta) h \|_{L^2 _{\mbox{{\scriptsize loc}}}}$. 
We expect the Remark~1 to extend to this case as well. 
The restriction $\gamma >0$ in Remark~1 seems more serious at first sight, 
since in the case $\gamma \le 0$, the coercivity estimate from Theorem~\ref{theo:Bolt} 
does not forbid the loss of mass at infinity. However on one hand for $\gamma =0$ 
and $b$ non locally integrable, the explicit diagonalization of the linearized 
Boltzmann operator (see~\cite{Bo88} for instance) shows that it has 
discrete spectrum (and compact resolvent) as well. On the other hand 
we shall give new improved coercivity estimates in weighted $L^2(M)$ spaces 
in the non-cutoff case in the forthcoming work~\cite{MoSt}, 
showing that this restriction can still be further relaxed.  
\smallskip

3. Our proof covers the physical case of inverse power-law 
interaction potentials $V(r)= \mbox{cst} \, r^{-(s-1)}$ with $s>5$ in dimension $N=3$ 
(for which $\gamma >0$ and $\alpha > 0$). The limit case $s=5$ corresponds 
to Maxwell molecules (for which $\gamma=0$ and $\alpha>0$). It 
can be treated thanks to the explicit diagonalization for Maxwell molecules 
(see the previous Remark~2). 
\smallskip 

\medskip

Concerning the linearized Landau operator we prove the 
 \begin{theorem}[The linearized Landau operator] \label{theo:Land}
 Under assumptions~\eqref{eq:hypPhi}, the linearized Landau operator 
 $L^{{\ron L}}$ with collision kernel $\Phi$ satisfies
   \begin{equation}\label{coLa}
   \forall \, h \in {\ron D}\left(L^{{\ron L}}\right), \quad  
   D^{{\ron L}}(h) \ge C^{{\ron L}} _\gamma \, 
             \left( \big\| \big[ h - \Pi(h)\big] \big\|_{H^1(\langle v \rangle ^\gamma M)} ^2 
                   + \big\| \big[ h - \Pi(h)\big] \, 
                                        \langle v \rangle ^{1+\gamma/2} \big\|_{L^2(M)} ^2 \right),
   \end{equation}   
 where $C^{{\ron L}} _\gamma$ is an explicit constant depending only on 
 $\gamma$, $C_\Phi$, and the dimension $N$. 
 \end{theorem}

\Remarks  

1. It was already noticed (by non-constructive arguments) 
in~\cite{Guo:Land} that the Dirichlet form of the linearized Landau operator 
controls the $L^2( \langle v \rangle ^{\gamma+2} M)$ norm. Moreover in~\cite{Guo:Land} it is given 
(still by non-constructive arguments) estimates of the form 
$$
D^{{\ron L}}(h) \ge C \, 
\left( 
\| h \|_{L^2(\langle v \rangle^{\gamma+2} M)}^2 + 
\left\| {\bf P}(v) \nabla_v h \right\|_{L^2(\langle v \rangle^{\gamma+2} M)} ^2 + 
\left\| (\mbox{Id} - {\bf P}(v)) \nabla_v h \right\|_{L^2(\langle v \rangle^\gamma M)} ^2
\right)
$$
which use on the right-hand side some norm slightly stronger than the one in~\eqref{coLa}. 
The proof of explicit coercivity estimates for this stronger norm shall be 
given in the forthcoming work~\cite{MoSt}.  
\smallskip

2. By a similar argument as for the linearized Boltzmann operator with 
a non locally integrable collision kernel, 
when $\gamma >-2$, we straightforwardly deduce 
from Theorem~\ref{theo:Land} that the linearized Landau operator 
has compact resolvent and thus a purely discrete spectrum. 
Indeed if one considers a sequence 
$(h_k)_{k \ge 0}$ which is bounded in $L^2(M)$ and such that 
$(L^{\cal{L}}(h_k))_{k \ge 0}$ is bounded in $L^2(M)$, the 
coercivity estimate~\eqref{coLa} implies that this sequence 
is bounded in $H^1 _{\mbox{{\scriptsize loc}}} \cap L^2(\langle v \rangle^{2 + \gamma} M)$. 
Since $2 +\gamma>0$ when $\gamma >-2$, this implies by Rellich-Kondrachov Theorem 
that the sequence has a cluster point in $L^2(M)$.  
\medskip

\subsection{Dependence of the constants in the coercivity estimates according to the equilibrium} \label{subsec:norm}

Straightforward computations show that the Dirichlet form with mass $\rho$, 
momentum $u$ and temperature $T$ satisfies (with obvious notation) 
$$ D^{{\ron B}} _{\rho,u,T} (h) = 
    \frac{ \rho^2 } { \pi^N} \, (2 \,T)^{\gamma/2} \, 
    D^{{\ron B}} \Big(h\big(u+\sqrt{2T}v\big)\Big) $$
in the Boltzmann case, and 
$$ D^{{\ron L}} _{\rho,u,T} (h) = 
    \frac{ \rho^2 } { \pi^N} \, (2 \,T)^{\gamma/2+1} \, 
    D^{{\ron L}} \Big(h\big(u+\sqrt{2T}v\big)\Big) $$
in the Landau case. Moreover we have for the norms involved in our estimates:  
$$ \left\| h \right\|_{L^2\big( \langle v \rangle ^\gamma M_{\rho,u,T}\big)} 
     =  \frac{ \rho } { \pi^{N/2}} \, (2 \,T)^{\gamma/2} \, \left\| h\big(u+\sqrt{2T}v\big) \right\|_{L^2( \langle v \rangle ^\gamma M)}$$ 
for the norms involved in the Boltzmann case (for the $H^{\alpha/2} _{\mbox{{\scriptsize loc}}}$ space 
the dependency cannot be written with such a simple formula) and 
$$ \left\| h \right\|_{H^1\big( \langle v \rangle ^\gamma M_{\rho,u,T}\big)} 
     =  \frac{ \rho } { \pi^{N/2}} \, (2 \,T)^{\gamma/2} \, \left\| h\big(u+\sqrt{2T}v\big) \right\|_{L^2( \langle v \rangle ^\gamma M)}$$ 
and 
$$ \left\| h \right\|_{L^2\big( \langle v \rangle ^{\gamma+2} M_{\rho,u,T}\big)} 
     =  \frac{ \rho } { \pi^{N/2}} \, (2 \,T)^{\gamma/2 +1} \, \left\| h\big(u+\sqrt{2T}v\big) \right\|_{L^2( \langle v \rangle ^\gamma M)}$$ 
for the norms involved in the Landau case. This explains how to modify the constants 
in Theorems~\ref{theo:Bolt} and~\ref{theo:Land} for a general equilibrium.  

For instance for an equilibrium with mass $\rho$, 
momentum $u$ and temperature $T$, the constant in Theorem~\ref{theo:Bolt} has to be 
multiplied by some factor $\rho/\pi^{N/2}$. This constant is independent of the momentum, which is 
a well-known consequence of the translation invariance of the Boltzmann equation. 
It is also independent of the temperature, which is due to the fact that our polynomial weight matches 
exactly the polynomial dependency of the collision kernel in terms of the relative velocity. 

\subsection{Method of proof}

In the case of hard potentials, the idea is to decompose the operator 
between a part satisfying the desired 
coercivity estimate and a bounded part, and use the spectral gap estimates. 
This argument is reminiscent of an argument of Grad \cite[Section~5]{Grad63} used 
to study the behavior at $v \to \infty$ of the eigenvectors of the linearized 
Boltzmann operator for hard potentials, and it was already noticed 
in~\cite{BCN:86}. 
Nevertheless it is the first 
time that it is used to obtain explicit estimates (thanks to the results in~\cite{BaMo}). 
The same idea, combined with a suitable Poincar\'e inequality, 
is applied to the linearized Landau operator. 


For soft potentials we decompose the Dirichlet form 
according to the modulus of the relative velocity. Combined with 
technical estimates on the non-local part of the linearized collision operators 
and the coercivity estimates from the Maxwell case, 
it enables to reconstruct a lower bound with the appropriate weight. 
The proof for the linearized Landau is strongly 
guided by the previous study of the Boltzmann case, which helps to 
identify relevant estimates. 

Finally the proof of the coercivity estimates in local Sobolev spaces for the linearized 
Boltzmann operator with a non locally integrable collision kernel 
is inspired by the previous works~\cite{Lions98,Vi:99,ADVW} 
on the full non-linear collision operator,  
and by our study of the linearized Landau operator. Indeed the 
suitable decomposition of $L^{{\ron B}}$ for non locally integrable 
collision kernels (for which the usual Grad's splitting does not 
make sense anymore) is directly readable on the 
linearized Landau operator: the part which becomes the diffusion part 
in the grazing collision limit is the part which enjoys a coercivity 
property in local Sobolev spaces, and the part which becomes the bounded 
part in the grazing collision limit is the part which is bounded thanks 
to the ``cancellation lemmas'' (which we borrow from~\cite{Vi:99,ADVW}). 
Let us also mention that, as indicated by one of the anonymous 
referee, a decomposition in the same spirit was proposed in 
the papers~\cite{Alex97,Alex99}. 

\subsection{Plan of the paper} 
Section~\ref{sec:Bolt} is devoted to the linearized Boltzmann operator: it 
contains the proof of Theorem~\ref{theo:Bolt}, 
divided into two parts, for hard and then soft potentials, and then 
the proof of Theorem~\ref{theo:Bolt:nc}. Section~\ref{sec:Land} is 
devoted to the linearized Landau operator: it contains the proof of 
Theorem~\ref{theo:Land}, divided into hard and soft potentials again. 

\section{The linearized Boltzmann operator}\label{sec:Bolt}
\setcounter{equation}{0}

In this section and the next one, the constants which are only internal to a proof shall be denoted 
$C_1,C_2,\dots$ if they are referred to inside the proof, or simply $C$ if not.  

\subsection{Hard potentials}

Notice that the case $\gamma=0$ of Theorem~\ref{theo:Bolt} is already proved by the 
explicit estimates of the Maxwell case, see \cite{Bo88}. 
Hence we assume that $\gamma >0$ and we pick $h \in L^2(M)$ orthogonal 
to the null space of $L^{{\ron B}}$. First using the minoration of $b$ \eqref{eq:hypb} we reduce to the 
(cutoff case) where $b \equiv 1$ by~\cite[Lemma~2.1]{BaMo}, and using the assumption~\eqref{eq:hypPhi} 
we reduce to the case $\Phi(z) = z^\gamma$.

We consider the decomposition 
  \[ L^{{\ron B}} = K^{{\ron B}} - A^{{\ron B}} \]
with 
  \[
  K^{{\ron B}} h(v)=\int_{\R^N \times \ens{S}^{N-1}} B(|v-v_*|,\cos \theta) \, 
  M(v_{*}) \left[ h^{'}_{*}+h^{'} - h_{*} \right] \, dv_{*} \, d\sigma
  \]
and 
  \[ 
  A^{{\ron B}} h(v)= \left( \int_{\R^N \times \ens{S}^{N-1}} B(|v-v_*|,\cos \theta) \, 
  M(v_{*}) \, dv_{*} \, d\sigma \right) \, h(v). 
  \]
Then we use Grad computations \cite[Sections 2, 3, 4]{Grad63} to obtain that 
$K^{{\ron B}}$ is a (compact) bounded operator (with explicit bound $C_K ^{{\ron B}}$) and $A^{{\ron B}}$ 
is the multiplication operator by the collision frequency $\nu$, given here by  
  \[ \nu (v) = |\ens{S}^{N-1}| \, \int_{\R^N} |v-v_*|^\gamma \, M(v_*) \, dv_*. \]
On one hand we have straightforwardly  
  \[ \int_{\R^N} (A^{{\ron B}}h) \, h \, M \, dv \ge C_1 \, \| h \, \langle v \rangle^{\gamma/2} \|^2 _{L^2(M)} \]
with $C_1>0$ depending on $\gamma$. 
On the other hand we know by~\cite[Theorem~1.1]{BaMo} 
that there is an explicit constant $C_2 >0$ such that 
  \[ D^{{\ron B}} (h) = - \int_{\R^N} (L^{{\ron B}} h) \, h \, M \, dv \ge C_2 \, \|h \|^2 _{L^2(M)}. \]
We deduce then that 
  \begin{multline*}
  \| h \, \langle v \rangle^{\gamma/2} \|^2 _{L^2(M)} 
      \le C_1 ^{-1} \, \int_{\R^N} (A^{{\ron B}}h) \, h \, M \, dv \\
  \le C_1 ^{-1} \, \left[ - \int_{\R^N} (L^{{\ron B}} h) \, h \, M \, dv + 
      \int_{\R^N} (K^{{\ron B}} h) \, h \, M \, dv \right] \\
  \le C_1 ^{-1} \, \left[ D^{{\ron B}} (h) + 
      C_K ^{{\ron B}} \, \|h\|_{L^2(M)} ^2 \right] \\
  \le C_1 ^{-1} \left[ 1 + C_K ^{{\ron B}} \, C_2 ^{-1} \right]
     \, D^{{\ron B}} (h) 
  \end{multline*}
which concludes the proof of Theorem~\ref{theo:Bolt} in the case $\gamma >0$. 

\subsection{Soft potentials}

We suppose now that $\gamma <0$ and we pick $h \in L^2(M)$ orthogonal 
to the null space of $L^{{\ron B}}$. First using \eqref{eq:hypb} we reduce to the 
(cutoff case) where $b \equiv 1$ by using~\cite[Lemma~2.1]{BaMo} again (this lemma is independent 
on the particular form of $\Phi$), and using \eqref{eq:hypPhi} 
we reduce to the case $\Phi(z) = \min \{z^\gamma, 1\}$.
\smallskip

\noindent
{\bf Step 1}. We need first a technical lemma on $K^{{\ron B}}$, in the case of Maxwell molecules.  
We define 
  \[ K_R ^{{\ron B}} h(v) = \int_{\R^N \times \ens{S}^{N-1}} {\bf 1}_{\{|v-v_*| \ge R\}} \, 
  M(v_{*}) \left[ h^{'}_{*}+h^{'} - h_{*} \right] \, dv_* \, d\sigma. \]
Then 
  \begin{lemma} \label{lemKB} 
  The bounded operator $K_R ^{{\ron B}}$ satisfies 
    \begin{equation*}
    ||| K_R ^{{\ron B}} |||_{L^2(M)} \xrightarrow[]{R \to \infty} 0 
    \end{equation*}
  with explicit rate (here $|||\cdot|||_{L^2(M)}$ denotes the usual operator norm on $L^2(M)$).
  \end{lemma}
\begin{proof}[Proof of Lemma~\ref{lemKB}]
First we decompose $K_R ^{{\ron B}} = T_R - U_R$ with (using the 
change of variable $\sigma \to - \sigma$ that exchanges $v'$ and 
$v' _*$)
  \[ T_R h (v) = 2 \int_{\R^N \times \ens{S}^{N-1}} {\bf 1}_{\{|v-v_*| \ge R\}} \, 
  M(v_{*}) h^{'} \, d\sigma \, dv_{*} \]
and 
  \[ U_R h(v) = \int_{\R^N \times \ens{S}^{N-1}} {\bf 1}_{\{|v-v_*| \ge R\}} \, 
  M(v_{*}) h_* \, d\sigma \, dv_{*}. \]

The proof for $U_R$ is straightforward:
  \[ \|U_R h \|_{L^2(M)} \le |\ens{S}^{N-1}|^{1/2} \, \|h\|_{L^2(M)} \, 
           \left( \int_{\R^N \times \R^N} {\bf 1}_{\{|v-v_*| \ge R\}} \, M \, M_* \, dv \, dv_* \right)^{1/2} \]
which gives the convergence to $0$ for the operator norm with the rate. 

The term $T_R$ is more tricky to handle. First, using $M M_* = M' M'_*$, we write it as 
  \[ T_R h (v) = 2 M(v)^{-1/2} \, \int_{\R^N \times \ens{S}^{N-1}} {\bf 1}_{\{|v-v_*| \ge R\}} \, 
  (M')^{1/2} h^{'} \, (M_*)^{1/2} (M' _*)^{1/2} \, d\sigma \, dv_{*}. \]
Then we use the bound
  \[ {\bf 1}_{\{|v-v_*| \ge R\}} \le {\bf 1}_{\{|v-v'| \ge R/\sqrt{2}\}} + {\bf 1}_{\{|v-v'_*| \ge R/\sqrt{2}\}} \]
which yields a corresponding decomposition $\big|T_R h\big| \le T^1 _R h + T^2 _R h$ with 
  \[ \left\{ 
     \begin{array}{l} \displaystyle 
     T^1 _R h (v) = 2 M^{-1/2} \, \int_{\R^N \times \ens{S}^{N-1}} {\bf 1}_{\{|v-v'| \ge R/\sqrt{2}\}} \, 
     (M')^{1/2} |h^{'}| \, (M_*)^{1/2} (M' _*)^{1/2} \, d\sigma \, dv_{*}, \vspace{0.2cm} \\ \displaystyle
     T^2 _R h (v) = 2 M^{-1/2} \, \int_{\R^N \times \ens{S}^{N-1}} {\bf 1}_{\{|v-v' _*| \ge R/\sqrt{2}\}} \, 
     (M')^{1/2} |h^{'}| \, (M_*)^{1/2} (M' _*)^{1/2} \, d\sigma \, dv_{*}. 
     \end{array}
     \right. \]   
Now we follow the computations by Grad \cite[Sections 2 and 3]{Grad63} 
(recalled in~\cite[Chapter 7, Section 2]{CIP94}) 
to compute and bound from above the kernel of these operators: we make the changes the variables 
\begin{itemize}
\item $\sigma \in \ens{S}^{N-1}$, $v_* \in \R^N$ $\longrightarrow$ $\omega = (v'-v)/|v'-v| \in \ens{S}^{N-1}$, 
$v_* \in \R^N$ (the jacobian is bounded by a constant); 
\item then $\omega \in \ens{S}^{N-1}$, $v_* \in \R^N$ $\longrightarrow$ $\omega \in \ens{S}^{N-1}$, $u=v-v_* \in \R^N$ 
(the jacobian is equal to $1$); 
\item then keeping $\omega$ fixed, decompose orthogonally $u = u_0 \omega +W$ with 
$u_0 \in \R$ and $W \in \omega^\bot$ (the jacobian is equal to $1$);  
\item finally keeping $W \in V^\bot$ fixed, $\omega \in \ens{S}^{N-1}$, $u_0 \in \R$ 
$\longrightarrow$ $V = u_0 \omega \in \R^N$ (the jacobian is $(1/2)|V|^{-(N-1)}$).
\end{itemize} 
We get thus
  \begin{multline*}
  \big|T^1 _R h(v)\big| \le C \, M(v)^{-1/2} \, 
        \int_{V \in \R^N} \int_{W \in V^\bot} \\ 
        |h(v+V)| M(v+V)^{1/2} \, {\bf 1}_{\{|V| \ge R/\sqrt{2}\}} \, 
          |V|^{-(N-1)} \,  M(v+W)^{1/2} \,  M(v+V+W)^{1/2} \, dV \, dW.
  \end{multline*}
Then using that ($\omega$ denotes the unit vector of $V$)
  \begin{multline*}
   M(v + W)^{1/2} \, M(v+V+W)^{1/2} = M(V)^{1/4} M(V+2(W+v))^{1/4} \\
      \le M(V)^{1/4} M(W+(v-(v\cdot \omega)\omega)), 
  \end{multline*}
we obtain the bound from above
  \[ \big|T^1 _R h(v)\big| \le C \, M(v)^{-1/2} \, 
        \int_{\R^N} |h(v+V)| M(v+V)^{1/2} \, {\bf 1}_{\{|V| \ge R/\sqrt{2}\}} \, 
          |V|^{-(N-1)} \, e^{-\frac{1}{4}|V|^2} \, dV. \]
By Young's inequality one deduces immediately the convergence to $0$ of $T^1 _R$ 
in the operator norm with explicit rate. 
On the other hand for $T^2 _R$ we use first that  
  \[ ||| {\bf 1}_{\{|\cdot| \ge r\}} T^2 _R |||_{L^2(M)} \le 
     ||| {\bf 1}_{\{|\cdot| \ge r\}} T _R |||_{L^2(M)} \le C \, (1+r)^{-1/2} \]
with explicit constant by Grad \cite[Section 4]{Grad63} (or see~\cite[Chapter 7, Section 2]{CIP94} again). 
Thus we pick $\var >0$ and then $r$ such that 
  \begin{equation}\label{T21}
  ||| {\bf 1}_{\{|\cdot| \ge r\}} T^2 _R |||_{L^2(M)} \le \var/2.
  \end{equation}
Then using again the changes of variables detailed above we get
  \begin{multline*}
  \big|{\bf 1}_{\{|v| \le r\}} T^2 _R h(v)\big| 
      \le C \, M(v)^{-1/2} \, {\bf 1}_{\{|v| \le r\}} \, 
        \int_{\R^N} |h(v+V)| M(v+V)^{1/2} \, |V|^{-(N-1)} \\ 
           \left[ \int_{V^\bot} M(v + W)^{1/2} \, M(v+V+W)^{1/2} \, {\bf 1}_{\{|W| \ge R/\sqrt{2}\}} \, dW \right] \, dV.
  \end{multline*}
We use that 
  \[ M(v + W)^{1/4} \, M(v+V+W)^{1/4} \le M(v + W)^{1/4} \le M(v)^{-1/4} M(W)^{1/8} \]
and 
  \[ M(v + W)^{1/4} \, M(v+V+W)^{1/4} \le M(V)^{1/8} \]
to obtain
  \begin{multline*}
  |{\bf 1}_{\{|v| \le r\}} T^2 _R h(v)| 
      \le C \, e^{3 r^2/4} \\
        \int_{\R^N} |h(v+V)| M(v+V)^{1/2} \, |V|^{-{N-1}} \, M(V)^{1/8} \, 
           \left[ \int_{V^\bot} M(W)^{1/8} \, {\bf 1}_{\{|W| \ge R/\sqrt{2}\}} \, dW \right] \, dV.
  \end{multline*}
Since the function $|V|^{-{N-1}} \, M(V)^{1/8}$ belongs to 
$L^1$, the convolution according to this function is bounded from $L^2$ into $L^2$, 
and we deduce that 
  \[ \| {\bf 1}_{\{|\cdot| \le r\}} T^2 _R \|_{L^2(M)} \le C_r \, 
        \left[ \int_{\R^{N-1}} M(W)^{1/8} \, {\bf 1}_{\{|W| \ge R/\sqrt{2}\}} \, dW \right] \, 
        \|h\|_{L^2(M)} \]
and thus, for $R$ big enough, 
  \[ ||| {\bf 1}_{\{|\cdot| \le r\}} T^2 _R |||_{L^2(M)} \le \var/2. \]
Together with~\eqref{T21} this shows that $T^2 _R$ goes to $0$ in the operator norm with 
explicit rate, which ends the proof. 
\end{proof}
\smallskip

\noindent
{\bf Step 2.} Let us do a dyadic decomposition of $D^{{\ron B}}(h)$. We fix 
a parameter $R>1$ and we use the following decomposition of identity:
  \[ {\bf 1} = {\bf 1}_{\{|u|\le R\}} + \sum_{n \ge 1} {\bf 1}_{\{R^n \le |u|\le R^{n+1}\}} \]
to obtain 
  \[ D^{{\ron B}} (h) \ge C \, \sum_{n \ge 0} R^{(n+1)\gamma} \, \tilde{D}^{{\ron B}}_n (h) \]
with 
  \[ \tilde{D}^{{\ron B}}_n (h) = \int_{\R^N \times \R^N \times \ens{S}^{N-1}} 
  {\bf 1}_{\{R^n \le |v-v_*|\le R^{n+1}\}} \, 
  \left[ h^{'}_{*}+h^{'}-h_{*}-h \right]^2 \, M \, M_{*}  \, dv \, dv_* \, d\sigma \]
for $n \ge 1$, and 
  \[ \tilde{D}^{{\ron B}}_0 (h) = \int_{\R^N \times \R^N \times \ens{S}^{N-1}} 
  {\bf 1}_{\{|v-v_*|\le R\}} \, 
  \left[ h^{'}_{*}+h^{'}-h_{*}-h \right]^2 \, M \, M_{*}  \, dv \, dv_* \, d\sigma. \]
Now if we define
  \[ D^{{\ron B}}_k (h) = \int_{\R^N \times \R^N \times \ens{S}^{N-1}} 
  {\bf 1}_{\{|v-v_*|\le R^{k+1}\}} \, 
  \left[ h^{'}_{*}+h^{'}-h_{*}-h \right]^2 \, M \, M_{*}  \, dv \, dv_* \, d\sigma \]
for any $k \ge 0$, we have 
  \begin{multline*}
  \sum_{k \ge 0} R^{(k+1)\gamma} \, D^{{\ron B}}_k (h) = 
      \sum_{k \ge 0} R^{(k+1)\gamma} \sum_{0 \le n \le k} \tilde{D}^{{\ron B}}_n (h) \\
  = \sum_{n \ge 0} \tilde{D}^{{\ron B}}_n (h) \left( \sum_{k \ge n} R^{(k+1)\gamma} \right) 
  = S \, \sum_{n \ge 0} R^{(n+1)\gamma} \tilde{D}^{{\ron B}}_n (h) 
  \end{multline*}
where the constant  
  \[ S = \sum_{k \ge 0} \big(R^\gamma \big)^k \]
is finite thanks to the fact that $R>1$. Thus we deduce that
  \[ D^{{\ron B}} (h) \ge \frac{C}{S} \, \sum_{n \ge 0} R^{(n+1)\gamma} \, D^{{\ron B}}_n (h). \]
\smallskip

\noindent
{\bf Step 3}. In this step we estimate each term of the dyadic decomposition. 
We fix $n_0 \in \N$ (to be latter chosen big enough) and we estimate $D^{{\ron B}}_n (h)$ 
for $n \ge n_0$. We denote $\chi_r$ the indicator function depending on the four variables 
$v,v_*,v',v_* '$ such that at least one of these four points belongs to $B(0,r)$. 
We also define the shorthand 
  \[ \Delta(F)=\big[ F' + F' _* - F_* - F \big]. \]
Then 
  \begin{multline*}
  D^{{\ron B}}_n (h) = \int_{\R^N \times \R^N \times \ens{S}^{N-1}} 
  {\bf 1}_{\{|v-v_*|\le R^{n+1}\}} \, 
  \Delta(h)^2 \, M \, M_{*}  \, dv \, dv_* \, d\sigma \\
  \ge \int_{\R^N \times \R^N \times \ens{S}^{N-1}} 
  {\bf 1}_{\{|v-v_*|\le R^{n+1}\}} \, \chi_r(v,v_*,v',v' _*) \, 
  \Delta(h)^2 \, M \, M_{*}  \, dv \, dv_* \, d\sigma. 
  \end{multline*}
We take $r=R^{n+2} - R^{n+1}$ and we denote $h_k = h \, {\bf 1}_{\{|\cdot|\le R^k\}}$. 
If one of the four collision points belongs to $B(0,R^{n+2} - R^{n+1})$ and the relative velocity 
is bounded by $R^{n+1}$, the collision sphere is included in 
$B(0,R^{n+2})$. Thus we deduce  
  \[
  D^{{\ron B}}_n (h) 
  \ge \int_{\R^N \times \R^N \times \ens{S}^{N-1}} 
  {\bf 1}_{\{|v-v_*|\le R^{n+1}\}} \, \chi_r(v,v_*,v',v' _*) \, 
  \Delta(h_{n+2})^2 \, M \, M_{*}  \, dv \, dv_* \, d\sigma. \]
Now we remove the indicator function $\chi_r$ by bounding from above 
the term corresponding to $1-\chi_r$, that is when all the four 
collision points have a modulus greater than $R^{n+2} - R^{n+1}$. 
Simple computations yield 
  \begin{multline*}
  D^{{\ron B}}_n (h)   \ge \int_{\R^N \times \R^N \times \ens{S}^{N-1}} 
  {\bf 1}_{\{|v-v_*|\le R^{n+1}\}} \, 
  \Delta(h_{n+2})^2 \, M \, M_{*}  \, dv \, dv_* \, d\sigma \\
  - C_1 \, e^{-(R^{n+2} - R^{n+1})} \, \|h_{n+2}\|^2 _{L^2(M)} 
  \end{multline*}
for an explicit constant $C_1 >0$. Then 
  \begin{multline*}
  \int_{\R^N \times \R^N \times \ens{S}^{N-1}} 
  {\bf 1}_{\{|v-v_*|\le R^{n+1}\}} \, 
  \Delta(h_{n+2})^2 \, M \, M_{*}  \, dv \, dv_* \, d\sigma \\
  = - 4 \, \int_{\R^N \times \R^N \times \ens{S}^{N-1}} 
  {\bf 1}_{\{|v-v_*|\le R^{n+1}\}} \, 
  h_{n+2} \, \big[ (h_{n+2})' +(h_{n+2})_* ' - (h_{n+2})_* \big] \, M \, M_{*} \, dv \, dv_* \, d\sigma \\
  + 4 \, \int_{\R^N \times \R^N \times \ens{S}^{N-1}} 
  {\bf 1}_{\{|v-v_*|\le R^{n+1}\}} \, 
  h_{n+2} ^2 \, M \, M_{*} \, dv \, dv_* \, d\sigma \\
  \ge - 4 \, \int_{\R^N \times \R^N \times \ens{S}^{N-1}} 
  h_{n+2} \, \big[ (h_{n+2})' +(h_{n+2})_* ' - (h_{n+2})_* \big] \, M \, M_{*} \, dv \, dv_* \, d\sigma \\
  - 4 \, \int_{\R^N} \left(K_{R^{n+1}} ^{{\ron B}} h_{n+2}\right) \, h_{n+2} \, M \, dv \\
  + 4 \, \int_{\R^N \times \R^N \times \ens{S}^{N-1}} 
  {\bf 1}_{\{|v-v_*|\le R^{n+1}\}} \, 
  h_n ^2 \, M \, M_{*} \, dv \, dv_* \, d\sigma,  
  \end{multline*}
and thus we deduce that 
  \begin{multline*}
  \int_{\R^N \times \R^N \times \ens{S}^{N-1}} 
  {\bf 1}_{\{|v-v_*|\le R^{n+1}\}} \, 
  \Delta(h_{n+2})^2 \, M \, M_{*} \, dv \, dv_* \, d\sigma \\
  \ge - 4 \, \int_{\R^N \times \R^N \times \ens{S}^{N-1}} 
  h_{n+2} \, \left[ (h_{n+2})' +(h_{n+2})_* ' - (h_{n+2})_* \right] \, M \, M_{*} \, dv \, dv_* \, d\sigma \\
  + 4 \, \int_{\R^N \times \R^N \times \ens{S}^{N-1}} 
  h_n ^2 \, M \, M_{*} \, d\sigma \, dv_{*} \, dv \\
  - 4 \, \int_{\R^N} \left(K_{R^{n+1}} ^{{\ron B}} h_{n+2}\right) \, h_{n+2} \, M \, dv  \\
  - 4 \, \int_{\R^N \times \R^N \times \ens{S}^{N-1}} 
  {\bf 1}_{\{|v-v_*|\ge R^{n+1}\}} \, 
  h_n ^2 \, M \, M_{*} \, dv \, dv_* \, d\sigma. 
  \end{multline*}
From Lemma~\ref{lemKB} we have 
  \[ - 4 \, \int_{\R^N} \left(K_{R^{n+1}} ^{{\ron B}} h_{n+2}\right) \, h_{n+2} \, M \, dv  
       \ge - \epsilon_1(R^{n+1}) \, \|h_{n+2}\|^2 _{L^2(M)} \]
where $\epsilon_1(r)$ is an explicit function going to $0$ as $r$ goes to infinity. 
Also when $v \in B(0,R^n)$ and $|v-v_*| \ge R^{n+1}$ we have by triangular inequality 
$|v_*| \ge R^{n+1} -R^n$, and thus simple computations show that 
  \[ 4 \, \int_{\R^N \times \R^N \times \ens{S}^{N-1}} 
  {\bf 1}_{\{|v-v_*|\ge R^{n+1}\}} \, 
  h_n ^2 \, M \, M_{*} \, dv \, dv_* \, d\sigma \le C_2 \, e^{-(R^{n+1} -R^n)} \, \|h_{n+2}\|^2 _{L^2(M)}. \]

Collecting every terms we deduce 
  \begin{multline*}
  \sum_{n \ge n_0} R^{(n+1)\gamma} \, D^{{\ron B}}_n (h) \ge \\
  \sum_{n \ge n_0} R^{(n+1)\gamma} \, \Bigg[ - 4 \, \int_{\R^N \times \R^N \times \ens{S}^{N-1}} 
  h_{n+2} \, \left[ (h_{n+2})' +(h_{n+2})_* ' - (h_{n+2})_* \right] \, M \, M_{*} \, dv \, dv_* \, d\sigma \\
  + 4 \, R^{2\gamma} \, \int_{\R^N \times \R^N \times \ens{S}^{N-1}} 
  h_{n+2} ^2 \, M \, M_{*} \, dv \, dv_*  \, d\sigma \\
  - C_1 \, e^{-(R^{n+2} - R^{n+1})} \, \|h_{n+2}\|^2 _{L^2(M)} - \epsilon_1(R^{n+1}) \, \|h_{n+2}\|^2 _{L^2(M)} \\
  - C_2 \, e^{-(R^{n+1} -R^n)} \, \|h_{n+2}\|^2 _{L^2(M)} \Bigg] 
  \end{multline*}
which writes 
  \begin{multline*}
  \sum_{n \ge n_0} R^{(n+1)\gamma} \, D^{{\ron B}}_n (h) \ge 
  \sum_{n \ge n_0}  R^{(n+1)\gamma} \, \Bigg[ \int_{\R^N \times \R^N \times \ens{S}^{N-1}} 
  \Delta(h_{n+2})^2 \, M \, M_{*} \, dv \, dv_* \, d\sigma \\
  - 4 \, (1-R^{2\gamma}) \, \| h_{n+2} \|^2 _{L^2(M)} 
  - C_1 \, e^{-(R^{n+2} - R^{n+1})} \, \|h_{n+2}\|^2 _{L^2(M)} \\ - \epsilon_1(R^{n+1}) \, \|h_{n+2}\|^2 _{L^2(M)} 
  - C_2 \, e^{-(R^{n+1} -R^n)} \, \|h_{n+2}\|^2 _{L^2(M)} \Bigg]. 
  \end{multline*}
Now we use the explicit spectral gap for Maxwell molecules to get 
  \[ \int_{\R^N \times \R^N \times \ens{S}^{N-1}} 
  \Delta(h_{n+2})^2 \, M \, M_{*} \, dv \, dv_* \, d\sigma 
  \ge \lambda \, \| h_{n+2} - \Pi(h_{n+2}) \|^2 _{L^2(M)} \]
for an explicit $\lambda >0$. Hence we deduce that 
  \begin{multline*}
  \sum_{n \ge n_0} R^{(n+1)\gamma} \, D^{{\ron B}}_n (h) \ge 
  \sum_{n \ge n_0} R^{(n+1)\gamma} \, \Bigg[ \lambda \, \|h_{n+2}\|^2 _{L^2(M)} 
  - \lambda \, \|\Pi(h_{n+2})\|^2 _{L^2(M)}  \\ - 4 \, (1-R^{2\gamma}) \, \| h_{n+2} \|^2 _{L^2(M)} 
  - C_1 \, e^{-(R^{n+2} - R^{n+1})} \, \|h_{n+2}\|^2 _{L^2(M)} \\
  - \epsilon_1(R^{n+1}) \, \|h_{n+2}\|^2 _{L^2(M)} 
  - C_2 \, e^{R^{n+1} -R^n} \, \|h_{n+2}\|^2 _{L^2(M)} \Bigg]. 
  \end{multline*}
Since $\Pi(h)=0$, we have 
  \begin{multline*}
  \|\Pi(h_{n+2})\|^2 = \big\| \Pi \big( h \,{\bf 1}_{\{|\cdot| \ge R^{n+2}\}}\big) \big\|^2 
     \le \left( \int_{\{|v| \ge R^{n+2}\}} |v|^{2+|\gamma|} \, M(v) \, dv \right) \, 
         \big\| h \langle v \rangle^{\gamma/2} \big\|^2 _{L^2(M)} \\
  \le C_3 \, e^{-R^{n+2}} \, \big\| h \langle v \rangle^{\gamma/2} \big\|^2 _{L^2(M)}. 
  \end{multline*}
Now if we choose $R-1>0$ small enough such that 
  \[ 4 \, (1-R^{2\gamma}) \le \frac{\lambda}{8}, \]
then $n_0$ big enough so that $R^{n+2} - R^{n+1}=R^{n+1} (R-1)$ and $R^{n+1} -R^n = R^n (R-1)$ big enough 
such that 
  \[  \forall \, n \ge n_0, \quad C_1 \, e^{-(R^{n+2} - R^{n+1})}, \ C_2 \, e^{-(R^{n+1} -R^n)} \ \le  \frac{\lambda}{8} \]
and also $n_0$ big enough such that $R^{n+1}$ big enough such that 
  \[ \forall \, n \ge n_0, \quad \epsilon_1(R^{n+1}) \le \frac{\lambda}{8}, \]
we obtain for this choice of $n_0$ and $R$:
  \begin{multline*}
  \sum_{n \ge n_0} R^{(n+1)\gamma} \, D^{{\ron B}}_n (h) \ge \frac{\lambda}{2} \,
  \sum_{n \ge n_0} R^{(n+1)\gamma} \,  \|h_{n+2}\|^2 _{L^2(M)}
  - C_3 \, \lambda \, \left( \sum_{n \ge n_0} e^{-R^{n+2}} \right) \, \| h \langle v \rangle^{\gamma/2} \|^2 _{L^2(M)}\\
  \ge \left[ C_4 \, R^{n_0 \gamma} - C_5 \, e^{-R^{n_0}} \right] \, \| h \langle v \rangle^{\gamma/2} \|^2 _{L^2(M)}
  \end{multline*} 
for some explicit constants $C_4,C_5 >0$ independent on $n_0$. Thus by taking $n_0$ large enough we deduce 
that 
  \[ \sum_{n \ge n_0} R^{(n+1)\gamma} \, D^{{\ron B}}_n (h) 
       \ge C_6 \, \big\| h \langle v \rangle^{\gamma/2} \big\|^2 \]
for some explicit constant $C_6>0$. Coming back to $D^{{\ron B}} (h)$, this ends 
the proof of Theorem \ref{theo:Bolt} in the case $\gamma <0$. 

\subsection{Regularity for long-range interactions}\label{subsec:nc}

We suppose here that the collision kernel $B$ satisfies~\eqref{eq:hyptens}, 
\eqref{eq:hypPhi}, \eqref{eq:hypbnc}, with $\alpha>0$ (the case $\alpha=0$ 
is deduced from Theorem~\ref{theo:Bolt}).   
Thus we are reduced to the case where $B(|v-v_*|,\cos \theta) = |v-v_*|^\gamma \, \theta^{-(N-1)-\alpha}$. 
By symmetrizing the Dirichlet form with the change of variable $\sigma \to -\sigma$, 
we can finally reduce to the case where 
$B(|v-v_*|,\cos \theta) = |v-v_*|^\gamma \, \theta^{-(N-1)-\alpha} \, {\bf 1}_{\theta \in [0,\pi/2]}$.  
We pick $h \in L^2(M)$ orthogonal to the null space of $L^{{\ron B}}$. 

We start by restricting the velocity variables to a bounded domain. Let us fix 
$R>0$, and let us denote by ${\cal I}_R$ a $C^\infty$ mollified indicator function of the variables 
$v,v_*$ which is $1$ on $B_R=B(0,R)$ and $0$ outside $B(0,R+1)$. 

We control from below the Dirichlet form by 
  \[ D^{{\ron B}}(h) \ge \frac{1}{4} \, \int_{\R^N \times \R^N \times \ens{S}^{N-1}} 
  B  \, {\cal I}_R \, \left[ h' + h' _*   
  - h -h_* \right]^2 \, 
  M\, M_{*} \, dv \, dv_* \, d\sigma \]
and we develop it as 
  \begin{multline}\label{eq:splitEB:nc}
  D^{{\ron B}} \ge \frac{1}{4} \, \int_{\R^N \times \R^N \times \ens{S}^{N-1}} 
  B \,{\cal I}_R \, \left( \big[ h'-h \big]^2  
  + \big[ h^{'}_{*}-h_* \big]^2 \right) \, 
  M\, M_{*} \, dv \, dv_* \, d\sigma \\
  + \frac{1}{2} \, \int_{\R^N \times \R^N \times \ens{S}^{N-1}} 
  B \, {\cal I}_R \,  \big( h'-h \big) 
  \big( h^{'}_{*}-h_* \big) \, 
  M\, M_{*} \, dv \, dv_* \, d\sigma.
  \end{multline}
The pre-postcollisional change of variable on the second term and 
the change of variable $(v,v_*,\sigma) \to (v_*,v,-\sigma)$ on the first term yield 
  \begin{multline*}
  D^{{\ron B}} \ge \frac{1}{2} \, \int_{\R^N \times \R^N \times \ens{S}^{N-1}} 
  B \,{\cal I}_R \, \big( h'-h \big)^2 \,  
  M\, M_{*} \, dv \, dv_* \, d\sigma \\
  - \int_{\R^N \times \R^N \times \ens{S}^{N-1}} 
  B \,{\cal I}_R \, h_* \, \big( h'-h \big) \, 
  M\, M_{*} \, dv \, dv_* \, d\sigma 
  =: I_1 ^R + I_2 ^R. 
  \end{multline*} 
Now we estimate separately $I_1 ^R$ from below and $I_2 ^R$ from above. 
For the term $I_1 ^R$, the Carleman representation (see \cite{Carl:57}) yields 
  \[ I_1 ^R \ge C \, \int_{B_R \times B_R} S(v,v') \, \frac{\big( h'-h \big)^2}{|v-v'|^{N+\alpha}} \, dv \, dv' \]
where
  \[ S(v,v') = M(v) \, 
      \int_{E_{v,v'} \cap B_R} {\bf 1}_{B_R}(v_*) 
      \, |v'-v'_*|^{1+\gamma+\alpha} \, {\bf 1}_{\{|v'-v| \le |v'_*-v|\}} \, M' _* \, dv'_* \]
and $E_{v,v'}$ is the hyperplan containing $v$ and orthogonal to $v-v'$ 
(for the derivation of this formula, see \cite[Section~4]{Vi:99}). 
The second indicator function in the formula for $S(v,v')$ comes from the 
restriction to $\theta \in [0,\pi/2]$ by the symmetrization above. 
It is easily seen that $S(v,v')$ is 
bounded from below by some constant $C>0$ on $B_R \times B_R$. It follows that 
  \[ I_1 ^R \ge C \, \int_{B_R \times B_R} \frac{\big( h'-h \big)^2}{|v-v'|^{N+\alpha}} \, dv \, dv' 
            \ge  C_1 \, \|h\|^2 _{H^{\alpha/2}(B_R)} \]
for some constant $C_1>0$ (for the last inequality see for instance~\cite{Adams}). 

As for the second term $I_2 ^R$, we use the change of variable 
of the cancellation lemma in~\cite[Section~3]{ADVW}: 
keeping $v_*$ fixed, change $v,\sigma$ into $v',\sigma$ (the jacobian is $\cos^{-N} \theta/2$). 
We obtain 
  \begin{multline*}
  I_2 ^R = \int_{\R^N \times \R^N \times \ens{S}^{N-1}} 
  B \, h_* \, h \, |v-v_*|^\gamma \, 
  M_* \\ 
  b(\cos \theta) \, \Big[ M(\psi_\sigma(v)) \, \cos^{-N-\gamma} \theta/2 \, {\cal I}_R(\psi_\sigma(v),v_*) - 
  M(v) \, {\cal I}_R(v,v_*) \Big]  dv \, dv_* \, d\sigma
  \end{multline*}
where $\psi_\sigma(v)$ is the transformation introduced in~\cite[Lemma~1]{ADVW}: 
it is the point in the plan defined by $v,v_*,\sigma$ such that 
  \[ (\psi_\sigma(v) - v) \, \bot \, (v_* -v) \ \mbox{ and } \  
     (\psi_\sigma(v) - v_*) \cdot (v-v_*) = \cos \theta/2. \]

Now let us fix $v$, $v_*$ and $\theta$. Then the modulus 
$|\psi_\sigma(v)-v|= \tan \theta/2 \, |v-v_*|$ is fixed, and the 
vector $\psi_\sigma(v)-v$ satisfies 
  \[ \psi_\sigma(v)-v = |\psi_\sigma(v)-v| \, \omega \]
where $\omega$ is the opposite of the unit vector supporting 
the projection of $\sigma$ on the plan orthogonal to $v-v_*$ 
(see \cite[Figure~1]{ADVW}). 
It motivates the study of quantities like  
 \begin{equation*}
 I(\varphi) = \int_{\ens{S}_{v-v_*} ^{N-2}} 
 \big(\varphi(v + \rho \omega)-\varphi(v) \big) \, d\omega
 \end{equation*}
where $\varphi$ denotes some $C^2$ function on $\R^N$, $\ens{S}_{v-v_*} ^{N-2}$ denotes 
the unit sphere in the plan orthogonal to $v-v_*$, and $\rho>0$.  
If $\nabla \varphi$ denotes the gradient of $\varphi$ and $\nabla^2 \varphi$ denotes its Hessian 
matrix, one has the following Taylor expansion: 
 \begin{equation*}
 \varphi(v + \rho \omega) = \varphi(v) + \rho \, \big( \nabla \varphi(v) \cdot \omega \big)
                     + \frac{\rho^2}{2} \, \left\langle \nabla^2 \varphi (v + \rho' \omega) 
                     \cdot \omega , \omega \right\rangle
 \end{equation*}
for some $0 \le \rho' \le \rho$. By bounding the last term and 
taking the integral over $\ens{S}_{v-v_*} ^{N-2}$, we get the estimate
 \begin{equation*}
 \left| I(\varphi) - \rho \, \left( \int_{\ens{S}_{v-v'_*} ^{N-2}} d\sigma \, 
     \big(\nabla \varphi(v) \cdot \sigma \big) \right) \right| 
 \le  \frac{\rho^2}{2} \, |\ens{S}^{N-2}| \, \|\varphi\|_{W^{2,\infty}}.
 \end{equation*}
As the term involving $\nabla \varphi$ vanishes by symmetry, we obtain 
 \begin{equation*}
 |I(\varphi)| \le \frac{\rho^2}{2} \, |\ens{S}^{N-2}| \, \|\varphi\|_{W^{2,\infty}}.
 \end{equation*}

We apply this computation to $\varphi(v)=M(v) \, {\cal I}(v,v_*)$ 
with $\rho=|\psi_\sigma(v)-v|= \tan \theta/2 \, |v-v_*|$ to find  
  \begin{multline*}
  \forall \, v,v_* \in B_R, \ \ \ 
        \left| \int_{\ens{S}^{N-1}} b(\cos \theta) \, \big[ M(\psi_\sigma(v)) {\cal I}(\psi_\sigma(v),v_*) 
          - M(v) {\cal I}(v,v_*) \big] \, d\sigma \right| \\ 
           \le C \, \int_0 ^{\pi/2} \theta^{-(N-1)-\alpha} \, \tan^2 \theta/2 \, \sin^{N-2} \theta/2 \, d\theta \le C_2
  \end{multline*}
for some finite constant $C_2 >0$. Finally we have immediately  
  \begin{multline*}
  \int_{\ens{S}^{N-1}} b(\cos \theta) \, \left| \cos^{-N-\gamma} \theta/2 - 1 \right| \, d\sigma \\
  \le \int_0 ^{\pi/2} \theta^{-(N-1)-\alpha} \, \left| \cos^{-N-\gamma} \theta/2 - 1 \right| \, \sin^{N-2} \theta/2 \, d\theta \le C_3 
  \end{multline*}
for some finite constant $C_3>0$. We thus deduce that 
  \[ \big|I_2 ^R\big| \le C \, \|h\|_{L^2(B_R)} ^2 \le C_4 \, \|h \langle v \rangle^{\gamma/2} \|^2 _{L^2(M)}. \]

Now we can conclude the proof of Theorem~\ref{theo:Bolt:nc}. For any $R>0$, we have 
  \[ \|h\|^2 _{H^{\alpha/2}(B_R)} \le C_1 ^{-1} \, I_1 ^R 
             \le C_1 ^{-1} \, \left[ D^{{\ron B}}(h) + |I^R _2| \right]
             \le C_1 ^{-1} \, \left[ D^{{\ron B}}(h) 
          + C_4 \, \|h \langle v \rangle^{\gamma/2} \|^2 _{L^2(M)} \right]. \]
Since $\Pi(h)=0$ we can use the coercivity estimate of Theorem~\ref{theo:Bolt} 
  \[ \|h \langle v \rangle^{\gamma/2} \|^2 _{L^2(M)} \le C_5 \, D^{{\ron B}}(h) \]
to deduce finally 
  \[ \|h\|^2 _{H^{\alpha/2}(B_R)} \le C_1 ^{-1} \, \left[ 1 + C_4 \, C_5 \right] \, D^{{\ron B}}(h). \]
Since it is valid for any $R>0$, we obtain 
  \[ D^{{\ron B}}(h) \ge C \, \|h\|^2 _{H^{\alpha/2} _{\mbox{{\scriptsize loc}}}} \]
for some explicit constant $C>0$. This ends the proof of Theorem \ref{theo:Bolt:nc}. 
\medskip

\Remark 
It is easy to see that in our proof the dependency on $R$ of 
the constant in the control of the $H^{\alpha/2} (B_R)$ norm 
by the Dirichlet form can be made explicit and grows exponentially 
in terms of $R$. More details on this point will be discussed in the forthcoming 
work~\cite{MoSt}. 
\medskip

\section{The linearized Landau operator}\label{sec:Land}
\setcounter{equation}{0}

Note that here on the contrary to~\cite{BaMo} we are not able to perform a 
grazing collision limit in the coercivity estimates for the linearized Boltzmann operator.  
Thus we do not try to deduce results on the linearized Landau operator from the 
Boltzmann case, instead we work directly on this operator. 
 
\subsection{Hard potentials and Maxwell molecules}

We consider $h \in L^2(M)$ orthogonal to the null space of $L^{{\ron L}}$,  
we assume that $\gamma \ge 0$ and, thanks to the assumption~\eqref{eq:hypPhi},  
we reduce to the case $\Phi(z)=z^\gamma$. 

Classical computations, which can be found in~\cite[Section~2]{DeLe97} for instance, show 
that the linearized Landau operator $L^{{\ron L}}$ decomposes as 
  \[ L^{{\ron L}} = K^{{\ron L}} - A^{{\ron L}} \]
where $K^{{\ron L}}$ is a (compact) bounded operator (with explicit bound $C_K ^{{\ron L}}$) 
and $A^{{\ron L}}$ is a diffusion operator whose Dirichlet form satisfies 
  \[ \int_{\R^N} (A^{{\ron L}}h) \, h \, M \, dv  
     = \int_{\R^N} \left(\nabla_v h\right)^t {\cal M}(v) 
         \left(\nabla_v h\right) \, M \, dv \]
where the matrix ${\cal M}$ is symmetric definite positive with its smallest eigenvalue bounded from 
below by $C \, \langle v \rangle^{\gamma}$ for an explicit 
constant $C >0$ (see~\cite[Section~2, Propositions 2.3 and 2.4]{DeLe97}). 
Thus we deduce that 
  \[ \int_{\R^N} ( A^{{\ron L}}h ) \, h \, M \, dv  \ge 
     C \, \int_{\R^N} |\nabla_v h|^2 \, \langle v \rangle^{\gamma} \, M \, dv. \]
First, we recall that, as noticed in~\cite{Lemo:00}, a simpler way to recover the coercivity result 
from~\cite[Section~3, Theorem 3.1]{DeLe97} is to apply the Bakry-Emery criterion (see~\cite[Chapter~9, 
Section~2]{Vi:tspt}), which implies that $M$ satisfies a Poincar\'e inequality with constant $2$, 
and thus (as $h$ has zero mean) 
  \[ \int_{\R^N} |\nabla_v h|^2 \, \langle v \rangle^{\gamma} \, M \, dv 
     \ge \int_{\R^N} |\nabla_v h|^2 \, M \, dv \ge 2 \, \int_{\R^N} h^2 \, M \, dv. \]
Now we want to obtain a stronger coercivity estimate. Thus we apply Bakry-Emery 
criterion to the measure 
  \[ m(v) = \langle v \rangle^{\gamma} \, M(v) = \exp \left[ - |v|^2 + \frac{\gamma}{2} \ln (1+|v|^2) \right] 
          =: \exp\left[ - \phi(v) \right]. \]
A straightforward computation shows that 
  \[ \nabla^2 \phi \ge (2-\gamma) \, Id \]
which implies, as $(2-\gamma) \ge 1$ thanks to the assumptions on $\gamma$, that 
$m$ satisfies a Poincar\'e inequality with constant $1$, and thus  
  \[  \int_{\R^N} |\nabla_v h|^2 \, \langle v \rangle^{\gamma} \, M \, dv 
      \ge \int_{\R^N} \left[h - 
          \left(\int_{\R^N} h \, \langle v \rangle^{\gamma} \, M \, dv\right) \right] ^2 
           \, \langle v \rangle^{\gamma} \, M \, dv. \]
Hence by developing 
  \[  \int_{\R^N} |\nabla_v h|^2 \, \langle v \rangle^{\gamma} \, M \, dv 
      \ge  \int_{\R^N} h^2 \, \langle v \rangle^{\gamma} \, M \, dv 
            -  \left(\int_{\R^N} h \, \langle v \rangle^{\gamma} \, M \, dv \right)^2. \]
Now as 
  \[ \left(\int_{\R^N} h \, \langle v \rangle^{\gamma} \, M \, dv \right)^2 \le C_1 \, 
              \left( \int_{\R^N} h^2 \, M \, dv \right) \]
for some explicit constant $C_1$, we deduce by collecting every term that 
  \[ \int_{\R^N} (A^{{\ron L}}h) \, h \, M \, dv  \ge C_2 \, 
                     \|h  \|_{H^1(\langle v \rangle^\gamma M)} ^2 
                    - C_3 \, \|h\|_{L^2(M)} ^2 \]
for some explicit constants $C_2,C_3 >0$. 

Besides we have by~\cite[Theorem~1.2]{BaMo} 
  \[ - \int_{\R^N} (L^{{\ron L}}h) \, h \, M \, dv \ge C_4 \, \|h\|_{L^2(M)} ^2 \]
for an explicit constant $C_4>0$. Now let us write
  \begin{multline*}
  \| h  \|^2 _{H^1(\langle v \rangle^\gamma M)} 
      \le C_2 ^{-1} \, \int_{\R^N} (A^{{\ron L}}h) \, h \, M \, dv 
          +  C_2 ^{-1} C_3 \, \|h\|_{L^2(M)} ^2 \\
  \le C_2 ^{-1} \, \left[ - \int_{\R^N} (L^{{\ron L}} h) \, h \, M \, dv + 
      \int_{\R^N} (K^{{\ron L}} h) \, h \, M \, dv + C_3 \, \|h\|_{L^2(M)} ^2\right] \\
  \le C_2 ^{-1} \, \left[ D^{{\ron L}} (h) + 
      (C_K ^{{\ron L}} + C_3) \, \|h\|_{L^2(M)} ^2 \right] 
  \le C_2 ^{-1} \, \left[ 1 + (C_K ^{{\ron L}} + C_3)\, C_4 ^{-1} \right] \, D^{{\ron L}} (h) 
  \end{multline*}
which concludes half of the proof of Theorem \ref{theo:Land} when $\gamma \ge 0$. 
It remains to control the $L^2(M)$ norm with weight $1+\gamma/2$. 

Let us denote $g = h M^{1/2}$. Then 
  \begin{multline*}
  \|h \|_{H^1(\langle v \rangle^\gamma M)} ^2 
    \ge \int_{\R^N} |\nabla h| ^2 M \, \langle v \rangle^\gamma \, dv = 
     \int_{\R^N} \left| \nabla \frac{g}{M^{1/2}} \right|^2 M \, \langle v \rangle^\gamma \, dv = \\
      \int_{\R^N} | \nabla g |^2 \, \langle v \rangle^\gamma \, dv 
      +  \int_{\R^N} |g|^2 \, \frac{|v|^2}4 \, \langle v \rangle^\gamma \, dv 
      + \int_{\R^N} g \nabla g \cdot v \, \langle v \rangle^\gamma \, dv \\
      = \int_{\R^N} | \nabla g |^2 \, \langle v \rangle^\gamma \, dv  
        +   \int_{\R^N} |g|^2 \, \frac{|v|^2}4 \, \langle v \rangle^\gamma \, dv 
        - \int_{\R^N} |g|^2 \, |\nabla ( v \langle v \rangle^\gamma)| \, dv  \\
      \ge \|h \, |v|^{1+\gamma/2} \|_{L^2 (M)} ^2 - C \, \|h \langle v \rangle^\gamma\|_{L^2 (M)} ^2 
   \end{multline*}
which implies immediately, combined to the previous estimate, 
inequality~\eqref{coLa} in Theorem \ref{theo:Land}.  

\subsection{Soft potentials}

We follow almost the same path as for the linearized Boltzmann operator. The 
starting point is the following coercivity estimate in the Maxwell case 
  \begin{multline}\label{SGLandmax}
  \frac{1}{2}\int_{\R^N \times \R^N} |v-v_*|^2 \, 
  \Big\|{\bf P}(v-v_*) \Big[ \left( \nabla h \right) - \left( \nabla h \right)_* \Big] \Big\|^2 \,  
  M\, M_{*} \, dv \, dv_* \\ 
  \ge \lambda \, \left[ \|h\|_{H^1(M)} ^2 + \|h \langle v \rangle \|_{L^2(M)} ^2 \right]
  \end{multline}
for some explicit constant $\lambda >0$, which has been proved in the previous subsection.  

We assume that $\gamma <0$ and we pick $h \in L^2(M)$ orthogonal 
to the null space of $L^{{\ron L}}$. Using the assumption~\eqref{eq:hypPhi} 
we reduce to the case $\Phi(z) = \min \{z^\gamma, 1\}$.
\smallskip

\noindent
{\bf Step 1}. We first prove a technical lemma on $K^{{\ron L}}$, in the case of Maxwell molecules.  
We define for $R>0$
  \[ K_R ^{{\ron L}} h(v) = - M(v) ^{-1} \, \nabla_v \cdot \left( \int_{\R^N} 
   \big(1-\Theta_R(v-v_*) \big) \, |v-v_*|^2 \, 
  {\bf P}(v-v_*) \left( \nabla h \right)_* \, M \, M_* \, dv_* \right) \]
where $\Theta_R$ is a $C^\infty$ function on $\R^N$ such that $0 \le \Theta_R \le 1$, 
$\Theta_R=1$ on $B(0,R)$ and $\Theta_R=0$ outside $B(0,R+1)$.  
Then
  \begin{lemma} \label{lemKBL} 
  The bounded operator $K_R ^{{\ron L}}$ satisfies 
    \begin{equation*}
    ||| K_R ^{{\ron L}} |||_{L^2(M)} \xrightarrow[]{R \to \infty} 0 
    \end{equation*}
  with explicit rate. 
  \end{lemma}
\begin{proof}[Proof of Lemma~\ref{lemKBL}]
It amounts to a differentiation under the integral, an integration by part, and 
Young's inequality. 
\end{proof}
\smallskip

\noindent
{\bf Step 2}. We do the same dyadic decomposition of $D^{{\ron L}}(h)$ as for 
the Boltzmann case, to obtain 
  \[ D^{{\ron L}} (h) \ge C \, \sum_{n \ge 0} R^{(n+1)\gamma} \, D^{{\ron L}}_n (h) \]
for some constant $C>0$, with 
  \[ D^{{\ron L}}_n (h) = \int_{\R^N \times \R^N} 
  {\bf 1}_{\{|v-v_*|\le R^{n+1}\}} \, |v-v_*|^2
  \big\|{\bf P}(v-v_*) \big[ \left( \nabla h \right) - \left( \nabla h \right)_*\big] \big\|^2 
  M\, M_{*} \, dv \, dv_* \]
for any $n \ge 0$. 
\smallskip

\noindent
{\bf Step 3}. In this step we estimate each term of the dyadic decomposition. 
We fix $n_0 \ge 0$ (to be later chosen big enough) and we estimate $D^{{\ron L}}_n (h)$ 
for $n \ge n_0$. We denote $\chi_r$ the indicator function depending on 
$v,v_*$ such that at least on of these two points belongs to $B(0,r)$. 
We also define the shorthand 
  \[ \Delta(F)=\big\|{\bf P}(v-v_*) \big[ \left( \nabla F \right) - \left( \nabla F \right)_*\big] \big\|. \]
Then 
  \begin{multline*}
  D^{{\ron L}}_n (h) = \int_{\R^N \times \R^N} 
  {\bf 1}_{\{|v-v_*|\le R^{n+1}\}} \, |v-v_*|^2 \, 
  \Delta(h)^2 \, M \, M_{*} \, dv \, dv_* \\
  \ge \int_{\R^N \times \R^N} 
  {\bf 1}_{\{|v-v_*|\le R^{n+1}\}} \, \chi_r(v,v_*) \,  |v-v_*|^2 \, 
  \Delta(h)^2 \, M \, M_{*} \, dv \, dv_*. 
  \end{multline*}
We take $r=R^{n+2} - R^{n+1}$ and we denote $h_k = h \, \bar \Theta^\eta _{R^k}$ 
where $\bar \Theta_{R} ^\eta$ is defined by  $0 \le \bar \Theta^\eta _R \le 1$, 
$\bar \Theta^\eta _R=1$ on $B(0,R)$, $\Theta_R=0$ outside $B(0,R+\eta^{-1})$ 
and $|\nabla \bar \Theta_R ^\eta| \le \eta$ (with $\eta>0$ to be later chosen small enough).  

If $v$ or $v_*$ belongs to $B(0,R^{n+2} - R^{n+1})$ and the relative velocity 
is bounded by $R^{n+1}$, both points belong to $B(0,R^{n+2})$. Thus we deduce 
  \[  D^{{\ron L}}_n (h) 
  \ge \int_{\R^N \times \R^N} 
  {\bf 1}_{\{|v-v_*|\le R^{n+1}\}} \, |v-v_*|^2 \, \chi_r(v,v_*) \, 
  \Delta(h_{n+2})^2 \, M \, M_{*} \, dv \, dv_*. \]
Now we remove the indicator function $\chi_r$ by bounding from above 
the term corresponding to $1-\chi_r$, that is when $v$ and $v_*$ 
have a modulus greater than $R^{n+2} - R^{n+1}$. 
Simple computations yield 
  \begin{multline*}
  D^{{\ron L}}_n (h)  \ge \int_{\R^N \times \R^N} 
  {\bf 1}_{\{|v-v_*|\le R^{n+1}\}} \, |v-v_*|^2 \, 
  \Delta(h_{n+2})^2 \, M \, M_{*} \, dv \, dv_* \\
  - C_1 \, R^{2(n+1)} \, e^{-(R^{n+2} - R^{n+1})} \, \|h_{n+2}\|^2 _{H^1(M)}
  \end{multline*}
for an explicit constant $C_1 >0$. Then we focus on the main term 
  \[ \int_{\R^N \times \R^N} 
  {\bf 1}_{\{|v-v_*|\le R^{n+1}\}} \, |v-v_*|^2 \,  
  \Delta(h_{n+2})^2 \, M \, M_{*} \, dv \, dv_*. \]
Since ${\bf 1}_{\{|v-v_*|\le R^{n+1}\}} \ge \Theta_{R^{n+1}-1}(v-v_*)$, we first bound it from below by   
  \[ \int_{\R^N \times \R^N} 
  \Theta_{R^{n+1}-1}(v-v_*) \, |v-v_*|^2 \,  
  \Delta(h_{n+2})^2 \, M \, M_{*} \, dv \, dv_*. \]

Then we proceed as in the Boltzmann case:  
  \begin{multline*}
  \int_{\R^N \times \R^N} 
  \Theta_{R^{n+1}-1}(v-v_*) \, |v-v_*|^2 \,  
  \Delta(h_{n+2})^2 \, M \, M_{*} \, dv \, dv_* \\ 
  = - 2 \, \int_{\R^N \times \R^N} 
  \Theta_{R^{n+1}-1}(v-v_*) \, |v-v_*|^2 \,  
  \big[ {\bf P}(v-v_*) \left( \nabla h_{n+2} \right)  \big]  \cdot  
  \big[ {\bf P}(v-v_*) \left( \nabla h_{n+2} \right)_*  \big] \, M \, M_{*} \, dv \, dv_* \\
  + 2 \, \int_{\R^N \times \R^N} 
  \Theta_{R^{n+1}-1}(v-v_*) \, |v-v_*|^2 \, 
  \big\|{\bf P}(v-v_*) \left( \nabla h_{n+2} \right)  \big\|^2 \, M \, M_{*} \, dv \, dv_* \\
  \ge - 2 \, \int_{\R^N \times \R^N} |v-v_*|^2 \, 
  \big[ {\bf P}(v-v_*) \left( \nabla h_{n+2} \right)  \big]  \cdot  
  \big[ {\bf P}(v-v_*) \left( \nabla h_{n+2} \right)_*  \big] \, M \, M_{*} \, dv \, dv_* \\
  - 2 \, \int_{\R^N} \left(K_{R^{n+1}-1} ^{{\ron L}} h_{n+2}\right) \, h_{n+2} \, M \, dv \\
  + 2 \, \int_{\R^N \times \R^N} \Theta_{R^{n+1}-1}(v-v_*) \, |v-v_*|^2 \, 
  \big\|{\bf P}(v-v_*) \left( \nabla h_{n+2} \right)  \big\|^2 \, M \, M_{*} \, dv \, dv_*.  
  \end{multline*}
Then we use that 
  \begin{multline*}
  \big\|{\bf P}(v-v_*) \left( \nabla h_{n+2} \right)  \big\|^2 
     =  \big\|{\bf P}(v-v_*) \left( \nabla (\bar \Theta_{R^{n+2}} ^\eta h) \right)  \big\|^2 \\
        \ge  (1-\var) \, \big\|{\bf P}(v-v_*) \left( \bar \Theta_{R^{n+2}} ^\eta \nabla h \right)  \big\|^2 
          - C(\var) \, \big\|{\bf P}(v-v_*) \left( h \nabla \bar \Theta_{R^{n+2}} ^\eta \right)  \big\|^2 \\ 
         \ge (1-\var) \, \big\|{\bf P}(v-v_*) \left( \bar \Theta_{R^n} ^\eta \nabla h \right)  \big\|^2 
          - C(\var) \, \big\|{\bf P}(v-v_*) \left( h \nabla \bar \Theta_{R^{n+2}} ^\eta \right)  \big\|^2  \\
         \ge (1-2\var) \, \big\|{\bf P}(v-v_*) \left( \nabla (\bar \Theta_{R^n} ^\eta h) \right)  \big\|^2 \\
          - C(\var) \, \big\|{\bf P}(v-v_*) \left( h \nabla \bar \Theta_{R^n} ^\eta \right)  \big\|^2 
          - C(\var) \,  \big\|{\bf P}(v-v_*) \left( h \nabla \bar \Theta_{R^{n+2}} ^\eta \right)  \big\|^2 \\ 
          \ge (1-2 \var) \, \big\|{\bf P}(v-v_*) \left( \nabla  h_n \right)  \big\|^2 
          - C_2(\var) \, \eta \, \big\| h_{n+3}  \big\|^2 
  \end{multline*}
since $R^{n+2} + \eta^{-1} \le R^{n+3}$ if $n \ge n_0$ is big enough. 

Hence we deduce that 
  \begin{multline*}
  \int_{\R^N \times \R^N} 
  \Theta_{R^{n+1}-1}(v-v_*)  \, |v-v_*|^2 \,  
  \Delta(h_{n+2})^2 \, M \, M_{*} \, dv_{*} \, dv \\
  \ge - 2 \, \int_{\R^N \times \R^N} |v-v_*|^2 \, 
  \left[ {\bf P}(v-v_*) \left( \nabla h_{n+2} \right)  \right]  \cdot  
  \left[ {\bf P}(v-v_*) \left( \nabla h_{n+2} \right)_*  \right] \, M \, M_{*} \, dv \, dv_*  \\
  + 2 \, (1-2\var) \, \int_{\R^N \times \R^N} |v-v_*|^2 \, 
  \left\|{\bf P}(v-v_*) \left( \nabla h_n \right)  \right\|^2 \, M \, M_{*} \, dv \, dv_* \\
  - 2 \, \int_{\R^N} \left(K_{R^{n+1}-1} ^{{\ron L}} h_{n+2}\right) \, h_{n+2} \, M \, dv 
  - C_2(\var) \, \eta \, \big\| h_{n+3} \langle v \rangle \big\|^2 _{L^2 (M)} \\ 
  - 2 \, \int_{\R^N \times \R^N} 
  {\bf 1}_{\{|v-v_*|\ge R^{n+1}-1\}} \,  |v-v_*|^2 \, 
  \left\|{\bf P}(v-v_*) \left( \nabla h_n \right)  \right\|^2 \, M \, M_{*} \, dv \, dv_*. 
  \end{multline*}
Now we use that (from Lemma~\ref{lemKBL})  
  \[  - 2 \, \int_{\R^N} \left(K_{R^{n+1}-1} ^{{\ron L}} h_{n+2}\right) \, h_{n+2} \, M \, dv  
       \ge - \epsilon_2(R^{n+1}-1) \, \|h_{n+2}\|^2 _{H^1(M)} \]
where $\epsilon_2(r)$ is an explicit function going to $0$ as $r$ goes to infinity. 
Also simple computations show that 
  \begin{multline*}
  2 \, \int_{\R^N \times \R^N} 
  {\bf 1}_{\{|v-v_*|\ge R^{n+1}-1\}} \,  |v-v_*|^2 \, 
  \left\|{\bf P}(v-v_*) \left( \nabla h_n \right)  \right\|^2 \, 
  M \, M_{*} \, dv \, dv_* \\ 
  \le C_3 \, R^{2n} \, e^{-(R^{n+1} -R^n)} \, \|h_{n+2}\|^2 _{H^1(M)}.
  \end{multline*}

Collecting every term we deduce 
  \begin{multline*}
  \sum_{n \ge n_0} R^{(n+1)\gamma} \, D^{{\ron L}}_n (h) \ge 
  \sum_{n \ge n_0} R^{(n+1)\gamma} \\ \Bigg[ - 2 \, \int_{\R^N \times \R^N} |v-v_*|^2 \, 
  \left[ {\bf P}(v-v_*) \left( \nabla h_{n+2} \right)  \right]  \cdot  
  \left[ {\bf P}(v-v_*) \left( \nabla h_{n+2} \right)_*  \right] \, M \, M_{*} \, dv \, dv_* \\
  + 2 \, (1-2\var)\, R^{2\gamma} \, \int_{\R^N \times \R^N} |v-v_*|^2 \, 
  \left\|{\bf P}(v-v_*) \left( \nabla h_{n+2} \right)  \right\|^2 \, M \, M_{*} \, dv \, dv_* \\ 
  - C_2(\var) \, \eta \, R^{-\gamma} \, \big\| h_{n+2}  \langle v \rangle \big\|^2 _{L^2 (M)} 
  - C_1 \, R^{2(n+1)} \, e^{-(R^{n+2} - R^{n+1})} \, \|h_{n+2}\|^2 \\ 
  - \epsilon_2(R^{n+1}-1) \, \|h_{n+2}\|^2 _{H^1(M)}
  - C_3 \, R^{2n} \, e^{-(R^{n+1} -R^n)} \, \|h_{n+2}\|^2 _{H^1(M)} \Bigg] 
  \end{multline*}
which writes 
  \begin{multline*}
  \sum_{n \ge n_0} R^{(n+1)\gamma} \, D^{{\ron L}}_n (h) \ge 
  \sum_{n \ge n_0} R^{(n+1)\gamma} \, \Bigg[ \int_{\R^N \times \R^N} |v-v_*|^2 \, 
  \Delta(h_{n+2})^2 \, M \, M_{*} \, dv \, dv_* \\ 
  - 2 \, (1-(1-2\var)R^{2\gamma}) \, \| h_{n+2} \|^2 _{H^1(M)} 
  - C_1 \, R^{2(n+1)} \, e^{-(R^{n+2} - R^{n+1})} \, \|h_{n+2}\|^2 _{H^1(M)} \\
  - C_2(\var) \, \eta \, R^{-\gamma} \, \big\| h_{n+2} \langle v \rangle \big\|^2 _{L^2 (M)} 
  - \epsilon_2(R^{n+1}) \, \|h_{n+2}\|^2 _{H^1(M)}
  - C_3 \, R^{2n} \, e^{-(R^{n+1} -R^n)} \, \|h_{n+2}\|^2 _{H^1(M)} \Bigg]. 
  \end{multline*}

Now we use the explicit coercivity estimate~\eqref{SGLandmax} for Maxwell molecules to deduce that 
  \begin{multline*}
  \sum_{n \ge n_0} R^{(n+1)\gamma} \, D^{{\ron L}}_n (h) \ge \\
  \sum_{n \ge n_0} R^{(n+1)\gamma} \, \Bigg[ \frac{3}{2} \, \lambda \, \|h_{n+2}\|^2 _{H^1(M)} 
   + \frac{3}{2} \, \lambda \, \|h_{n+2} \langle v \rangle \|^2 _{L^2(M)} 
  - 6 \, \lambda \, \|\Pi(h_{n+2})\|^2 _{H^1(M)} \\ 
  - 6 \, \lambda \, \|\Pi(h_{n+2}) \langle v \rangle \|^2 _{L^2(M)} 
  - 2 \, (1-(1-2\var)R^{2\gamma}) \, \| h_{n+2} \|^2 _{H^1(M)} \\
  - C_2(\var) \, \eta \, R^{-\gamma} \, \big\| h_{n+2} \langle v \rangle \big\|^2 _{L^2 (M)}
  - C_1 \, R^{2(n+1)} \, e^{-(R^{n+2} - R^{n+1})} \, \|h_{n+2}\|^2 _{H^1(M)} \\
  - \epsilon_2(R^{n+1}) \, \|h_{n+2}\|^2 _{H^1(M)}
  - C_3 \, R^{2n} \, e^{-(R^{n+1} -R^n)} \, \|h_{n+2}\|^2 _{H^1(M)} \Bigg]. 
  \end{multline*}
Since $\Pi(h)=0$, we have 
  \begin{multline*}
  6 \, \lambda \, \|\Pi(h_{n+2})\|^2 _{H^1(M)} 
  + 6 \, \lambda \, \|\Pi(h_{n+2}) \langle v \rangle \|^2 _{L^2 (M)} \\
      = \big\| \Pi ( h \,{\bf 1}_{\{|\cdot| \ge R^{n+2}\}}) \big\|^2 _{H^1(M)} 
     + \big\| \Pi ( h \,{\bf 1}_{\{|\cdot| \ge R^{n+2}\}}) \langle v \rangle \big\|^2 _{L^2 (M)} \\
  \le C_4 \, e^{-R^{n+2}} \, \big\| h \big\|^2 _{H^1(\langle v \rangle^\gamma M)}.
  \end{multline*}

Now if we choose $\var$ and $R-1>0$ small enough such that 
  \[ 2 \, (1-(1-2\var)R^{2\gamma}) \le \frac{\lambda}{4},  \]
then $\eta$ small enough such that 
  \[ C_2(\var) \, \eta \le  \frac{\lambda}{4}, \]
then $n_0$ big enough so that $R^{n+2} (R-1) \ge \eta^{-1}$ for any $n \ge n_0$ 
(see the discussion above), and so that 
$R^{n+2} - R^{n+1}=R^{n+1} (R-1)$ and $R^{n+1} -R^n = R^n (R-1)$ big enough such that 
  \[ \forall \, n \ge n_0, \quad C_1 \, R^{2(n+1)} \, e^{-(R^{n+2} - R^{n+1})}, \ 
     C_3 \, R^{2n} \, e^{-(R^{n+1} -R^n)} \ \le  \frac{\lambda}{4},  \]
and also $n_0$ big enough such that $R^{n+1}$ big enough such that 
  \[ \forall \, n \ge n_0, \quad \epsilon_2(R^{n+1}) \le \frac{\lambda}{4}, \]
we obtain for this choice of $R$, $\eta$ and $n_0$:
  \begin{multline*}
  \sum_{n \ge n_0} R^{(n+1)\gamma} \, D^{{\ron L}}_n (h) \\ \ge \frac{\lambda}{4} \,
  \sum_{n \ge n_0} R^{(n+1)\gamma} \, \left(  \|h_{n+2}\|^2 _{H^1(M)} 
     +  \|h_{n+2}\|^2 _{L^2( \langle v \rangle^2 M)} \right) \\ 
  - C_4 \, \lambda \, \left( \sum_{n \ge n_0} e^{-R^{n+2}} \right) \, 
  \big\| h \big\|^2 _{H^1(\langle v \rangle^\gamma M)} \\
  \ge \left[ C_5 \, R^{n_0 \gamma} - C_6 \, e^{-R^{n_0}} \right] \, 
  \left( \big\| h \big\|^2 _{H^1(\langle v \rangle^\gamma M)} + 
         \big\| h \langle v \rangle^{1+ \gamma/2} \big\|^2 _{L^2(M)} \right) 
  \end{multline*} 
for some explicit constants $C_5,C_6 >0$ independent on $n_0$. 
Thus by taking $n_0$ large enough we deduce that 
  \[ \sum_{n \ge n_0} R^{(n+1)\gamma} \, D^{{\ron L}}_n (h) \ge 
              C_7 \, \left( 
     \big\| h \big\|^2 _{H^1(\langle v \rangle^\gamma M)} 
     + \big\| h \langle v \rangle^{1+ \gamma/2} \big\|^2 _{L^2(M)} \right)  \]
for some explicit constant $C_7>0$. Coming back to $D^{{\ron L}} (h)$, this concludes the 
proof of Theorem \ref{theo:Land} when $\gamma <0$. 
\medskip

\bigskip
\noindent
{\bf{Acknowledgment}}: We thank Fran\c cois Golse for pointing us 
reference~\cite{BCN:86}, and we thank Yan Guo and Robert Strain for useful 
discussions on the Landau operator, and pointing us the results 
in~\cite{Guo:Land}. We also thank the three anonymous referees for their  
numerous comments and suggestions. Support by the European network HYKE, 
funded by the EC as contract HPRN-CT-2002-00282, is acknowledged. 

\begin{flushright} \signcm \end{flushright}

\end{document}